\RequirePackage{fix-cm}
\documentclass[smallextended]{svjour3}       % onecolumn (second format)
\smartqed  % flush right qed marks, e.g. at end of proof
\usepackage{graphicx}
\usepackage{amsmath}
\usepackage{amssymb}
\usepackage{latexsym}
\usepackage{mathrsfs}
\usepackage{colortbl}
\usepackage{amscd}
\usepackage{tikz-cd}
\usepackage{comment}

\journalname{Calcolo}
\def\diam{\mathop{\mathrm{diam}}}
\def\diag{\mathop{\mathrm{diag}}}
\def\>{\textgreater}
\def\<{\textless}

\def\conv{\mathop{\mathrm{conv}}}

\spnewtheorem{thr}{Theorem}{\bf}{\it}
\spnewtheorem{coro}{Corollary}{\bf}{\it}
\spnewtheorem{defi}{Definition}{\bf}{\it}
\spnewtheorem{lem}{Lemma}{\bf}{\it}
\spnewtheorem{prop}{Proposition}{\bf}{\it}
\spnewtheorem{ass}{Assumption}{\bf}{\it}
\spnewtheorem{cond}{Condition}{\bf}{\it}
\spnewtheorem{Rem}{Remark}{\it}{\rm}
\spnewtheorem{Ex}{Example}{\it}{\it}
\spnewtheorem{Note}{Note}{\it}{\it}
\spnewtheorem*{ex*}{Example:}{\it}{\it}
\spnewtheorem*{pf*}{Proof}{\bf}{\rm}
\spnewtheorem*{rem*}{Remark:}{\it}{\it}
\spnewtheorem*{note*}{Note:}{\it}{\it}

\spnewtheorem*{lem1*}{Lemma 1}{\bf}{\rm}
\spnewtheorem*{lem3*}{Lemma 3}{\bf}{\rm}

\spnewtheorem*{th21*}{Theorem A}{\bf}{\it}
\spnewtheorem*{thrBA*}{Theorem B}{\bf}{\it}
\spnewtheorem*{thrF*}{Theorem D}{\bf}{\it}
\spnewtheorem*{thrin*}{Theorem C}{\bf}{\it}
\spnewtheorem*{thrlp*}{Theorem E}{\bf}{\it}

\spnewtheorem*{th2*}{Theorem 2}{\bf}{\it}

\spnewtheorem*{coex1*}{Counterexample 1}{\it}{\it}
\spnewtheorem*{coex2*}{Counterexample 2}{\it}{\it}

\allowdisplaybreaks[3]

\newcounter{sone}
\newcounter{stwo}
\newcounter{sthree}
\newcounter{sfour}
\newcounter{sfive}
\newcounter{ssix}
\newcounter{lone}
\newcounter{ltwo}
\newcounter{lthree}
\newcounter{lfour}
\newcounter{lfive}
\newcounter{lsix}
\makeatletter
    
    \@addtoreset{equation}{section}
  \makeatother

\begin{document}

%\title{General theory of interpolation error estimates on anisotropic meshes, part \Roman{lthree}: The Raviart--Thomas finite element}%\thanks{Grants or other notes
%placed here. General acknowledgments should be placed at the end of the article.}}
%\title{Remarks on interpolation error analysis and inverse inequalities on anisotropic meshes}
%\title{Remarks on anisotropic interpolation error analysis}
%\title{Classification of interpolation error estimates on anisotropic meshes}
%\title{The Raviart--Thomas and the Morley interpolation error estimates on anisotropic meshes}
\title{Anisotropic Raviart--Thomas interpolation error estimates using a new geometric parameter}

\titlerunning{Interpolation error analysis on anisotropic meshes}        % if too long for running head

\author{Hiroki Ishizaka 
%\author{Hiroki Ishizaka \and Kenta Kobayashi \and Takuya Tsuchiya %etc.
}


\institute{Hiroki Ishizaka \at
              Graduate School of Science and Engineering, Ehime University, Matsuyama, Japan \\
              \email{h.ishizaka005@gmail.com}           %  \\
           \and
          % Kenta Kobayashi\at
             % Graduate School of Business Administration, Hitotsubashi University, Kunitachi, Japan \\
            %\email{kenta.k@r.hit-u.ac.jp}
            \and
            %Takuya Tsuchiya \at 
             %Graduate School of Science and Engineering, Ehime University, Matsuyama, Japan \\
              %\email{tsuchiya@math.sci.ehime-u.ac.jp}  
}

\date{Received: date / Accepted: date}

\maketitle

\begin{abstract}
We present precise Raviart--Thomas interpolation error estimates on anisotropic meshes. The novel aspect of our theory is the introduction of a new geometric parameter of simplices. It is possible to obtain new anisotropic Raviart--Thoma error estimates using the parameter. We also include corrections to an error in ``General theory of interpolation error estimates on anisotropic meshes" (Japan Journal of Industrial and Applied Mathematics, 38 (2021) 163-191), in which Theorem 3 was incorrect.  

%The feature of our theory is that the errors are bound in terms of the diameter and a new geometric parameter of a simplex. 
%This paper shows that stability and error estimates of the global Raviart--Thomas interpolation are obtained if the ratio of these quantities is bounded.
%The boundedness of the ratio of these quantities is equivalent to the maximum-angle conditions. 
%Through the introduction of the geometric parameter, stability and error estimates of the global Raviart--Thomas interpolation is obtained if the maximum-angle conditions hold. 
%This paper also includes corrections to an error in "General theory of interpolation error estimates on anisotropic meshes" (Japan Journal of Industrial and Applied Mathematics, 38 (2021) 163-191), in which Theorem 3 was incorrect.  

\end{abstract}

\section{Introduction}
%Interpolation error analysis is an important theme of finite element analysis. It is well-known that the shape-regular condition for mesh partitions leads to the optimal interpolation error estimates, which gives the convergence of finite element approximation for partial differential equations. However, the geometric condition can be relaxed to the \textit{maximum-angle condition} which makes it possible to use anisotropic meshes, see \cite{BabAzi76} for two-dimensional cases and \cite{Kri92} for three-dimensional cases. The use of anisotropic meshes is effective for, for example, problems in which the solution has anisotropic behaviour in some direction of the domain.
Interpolation error analysis is an essential theme in finite element analysis. It is well known that the shape-regular condition for mesh partitions leads to optimal interpolation error estimates, which results in the convergence of finite element approximation for partial differential equations. However, the geometric condition can be relaxed to the \textit{maximum-angle condition}, which makes it possible to use anisotropic meshes; see \cite{BabAzi76} for two-dimensional cases and \cite{Kri92} for three-dimensional cases. Anisotropic meshes are effective for, for example, problems in which the solution has anisotropic behaviour in some direction of the domain.

%This paper treats the Raviart--Thomas interpolation error estimates on anisotropic meshes. The Raviart--Thomas finite element space was proposed in \cite{RavTho77}. The space is used for mixed finite element methods of second order elliptic problems and incompressible flow problems. Anisotropic Raviart--Thomas interpolation theory has been developed \cite{AcoDur99,AcoApe10}. Meanwhile, in a recent paper \cite{IshKobTsu21a}, the authors proposed the new geometric parameter $H_{T}$ of $d$-simplices (see Definition \ref{defiH} for the definition) and derived Raviart--Thomas interpolation error estimates:
%In contrast, the authors derived Raviart--Thomas interpolations on $d$-simplices, and proposed a new parameter $H_{T}$ (see Definition \ref{defiH} for the definition) in a recent paper \cite{IshKobTsu21a}. The new parameter was used to derive
This paper considers Raviart--Thomas interpolation error estimates on anisotropic meshes. The Raviart--Thomas finite element space was proposed in \cite{RavTho77}. The space is used for mixed finite element methods for second-order elliptic and incompressible flow problems. Anisotropic Raviart--Thomas interpolation theory was developed \cite{AcoDur99,AcoApe10}. The key idea is to derive the component-wise stabilities of the  Raviart--Thomas interpolation on reference elements (Lemmata \ref{rt=lem7} and \ref{rt=lem9}).

In contrast, this paper proposes anisotropic Raviart--Thomas interpolation error estimates using a new parameter $H_{T_0}$ of $d$-simplices, $d \in \{2,3 \}$, (see Definition \ref{defi1}) proposed in a recent paper \cite{IshKobTsu21a} under an assumption and using the component-wise stabilities, for example, we derive the following anisotropic error  estimate (Theorem \ref{thr1161}):
\begin{align*}
\displaystyle
&\| I_{T_0}^{RT^k} v_0 - v_0 \|_{L^p(T_0)^d} \notag \\
&\leq c  \Biggl( \frac{H_{T_0}}{h_{T_0}}  \sum_{|\varepsilon| = \ell + 1} \mathscr{H}^{\varepsilon} \| \partial^{\varepsilon}_{x} (\Psi_{T_0}^{-1} v_0) \|_{L^p(\Phi_{T_1}^{-1}(T_0))^d} \notag \\
&\hspace{1cm} +  h_{T_0} \sum_{|\beta| = \ell} \mathscr{H}^{\beta} \| \partial^{\beta}_{x} \nabla_{x} \cdot (\Psi_{T_0}^{-1} v_0) \|_{L^{p}(\Phi_{T_0}^{-1}(T_0))} \Biggr),
\end{align*}
also see Theorem \ref{thr1161b}, where $I_{T_0}^{RT^k}:  W^{1,1}(T_0)^d \to RT^k(T_0)$ is the Raviart--Thomas interpolation on $T_0$ defined by \eqref{RTinter1} and \eqref{RTinter2}, $\Phi_{T_0}$ and $\Psi_{T_0}$ are respectively the Affine mapping and Piola transformation defined in Section \ref{Affinedef} and Section \ref{Pioladef}, $\varepsilon := (\varepsilon_1 , \ldots, \varepsilon_d) \in \mathbb{N}_0^d$ ($\mathbb{N}_0 := \mathbb{N} \cup \{ 0 \}$) is a multi-index, and $\mathscr{H}^{\varepsilon}$ is defined in Section \ref{addinot}. We remark that the heart of the proof of Theorems \ref{thr1161} and \ref{thr1161b} is the \textit{scaling argument} described in Section \ref{scaling=sec}. We are naturally able to consider the following geometric condition, which is equivalent to the maximum-angle condition, as being sufficient to obtain optimal order estimates: there exists $\gamma_0 \> 0$ such that
\begin{align}
\displaystyle
\frac{H_{T_0}}{h_{T_0}} \leq \gamma_0 \quad \forall \mathbb{T}_h \in \{ \mathbb{T}_h \}, \quad \forall T_0 \in \mathbb{T}_h, \label{NewGeo}
\end{align}
where $\mathbb{T}_h = \{ T_0 \}$ is a simplicial mesh of a domain in $\mathbb{R}^d$. The new geometric condition appears to be simpler than the maximum-angle condition. Furthermore, the quantity $\frac{H_{T_0}}{h_{T_0}}$ can be easily calculated in the numerical process of finite element methods. Therefore, the new condition may be useful. We expect that the Raviart--Thomas interpolation error estimates that include the new mesh condition will be effective for a posteriori error estimates. We call the element $T_0$ "good" when there exists a constant $\gamma_0 \> 0$ that satisfies the new geometric condition (see Remark \ref{goodel}). Arguments about "good elements" can be found in \cite{phd,IshKobTsu21c}.

We describe the previous paper \cite{IshKobTsu21a,IshKobTsu21b}. The authors derived Raviart--Thomas interpolations on $d$-simplices. However, the proof of Theorem 3 in \cite{IshKobTsu21a} included a mistake, and we need to modify its statement to correct this error (Section \ref{Sec62}). The Babu\v{s}ka and Aziz technique is generally not applicable on anisotropic meshes in the proof of Theorem 3 in \cite{IshKobTsu21a}. Although we argued in \cite{IshKobTsu21a,IshKobTsu21b} that we do not impose either the shape-regular or maximum-angle condition during mesh partitioning, we realised that the new geometric condition \eqref{NewGeo} is necessary for the stability results of the global Raviart--Thomas interpolation.

When there is no ambiguity, we use the notation and definitions given in \cite{IshKobTsu21a}. Throughout this paper, $c$ denotes a constant independent of $h$ (defined later) unless specified otherwise. These values may change in each context. Let $\mathbb{R}_+$ be the set of positive real numbers. For $k \in \mathbb{N}_0$, $\mathcal{P}^k(T)$ is spanned by the restriction to $T$ of polynomials in $\mathcal{P}^k$, where  $\mathcal{P}^k$ denotes the space of polynomials with the degree at most $k$. Furthermore, we often use the following inequality (see \cite[Exercise 1.20]{ErnGue04}). Let $0 \< r \leq s$ and $a_i \geq 0$, $i=1,2,\ldots,n$ ($n \in \mathbb{N}$), be real numbers. We then have
\begin{align}
\displaystyle
\left( \sum_{i=1}^n a_i^s \right)^{1/s} \leq \left( \sum_{i=1}^n a_i^r \right)^{1/r}. \label{jensen}
\end{align}	
For matrix $A \in \mathbb{R}^{d \times d}$, we denote by $[{A}]_{ij}$ the $(i,j)$-component of $A$. We set $\| A \|_{\max} := \max_{1 \leq i,j \leq d} | [{A}]_{ij} |$. We also use the inequality
\begin{align}
\displaystyle
 \| {A} \|_{\max} \leq \| {A} \|_2. \label{Anorm}
 \end{align}

\section{Settings for the analysis of anisotropic interpolation theory}
Throughout this paper, let $d \in \{ 2 , 3 \}$. Let $\Omega \subset \mathbb{R}^d$ be a bounded polyhedral domain. Let $\mathbb{T}_h = \{ T_0 \}$ be a simplicial mesh of $\overline{\Omega}$ made up of closed $d$-simplices, such as
\begin{align*}
\displaystyle
\overline{\Omega} = \bigcup_{T_0 \in \mathbb{T}_h} T_0,
\end{align*}
with $h := \max_{T_0 \in \mathbb{T}_h} h_{T_0}$, where $ h_{T_0} := \diam(T_0)$. For simplicity, we assume that $\mathbb{T}_h$ is conformal: that is, $\mathbb{T}_h$ is a simplicial mesh of $\overline{\Omega}$ without hanging nodes.

\subsection{Motivation}
Let  $T_0 \subset \mathbb{R}^d$ and $\widehat{T} \subset \mathbb{R}^d$  be an element on a mesh $\mathbb{T}_h$ and a reference element defined in Section \ref{reference}. Let these two elements be affine equivalent. As an usual manner, the transformation $\Phi_0$ takes the form
\begin{align*}
\displaystyle
&\Phi_0: \widehat{T} \ni \hat{x}  \mapsto \Phi_0(\hat{x}) := {B}_0 \hat{x} + b_{0} \in T_0, 
\end{align*}
where ${B}_{0} \in \mathbb{R}^{d \times d}$ is an invertible matrix and $b_{0} \in \mathbb{R}^{d}$. According to the classical theory (e.g., see \cite[Theorem 1.114]{ErnGue04}), it holds that
\begin{align*}
\displaystyle
\| v - I_{T_0}^{RT^0} v \|_{L^{p}(T_0)^d} &\leq c \left( \| {B}_{0} \|_2 \| {B}_{0}^{-1} \|_2 \right ) \| {B}_{0} \|_2 | v |_{W^{1,p}(T_0)^d} \quad \forall v \in W^{1,p}(T_0)^d.
\end{align*}
Here, the quantity $ \| {B}_{0} \|_2 \| {B}_{0}^{-1} \|_2$ is called the \textit{Euclidean condition number} of ${B}_{0}$. By standard estimates (e.g., see \cite[Lemma 1.100]{ErnGue04}), we have
\begin{align*}
\displaystyle
\| {B}_{0} \|_2 \| {B}_{0}^{-1} \|_2 \leq c \frac{h_{T_0}}{\rho_{T_0}}, \quad  \| {B}_{0} \|_2 \leq c h_{T_0},
\end{align*}
where $\rho_{T_0}$ is the diameter of the largest ball that can be inscribed in $T_0$. It thus holds that
\begin{align*}
\displaystyle
\| v - I_{T_0}^{RT^0} v \|_{L^{p}(T_0)^d} &\leq c \frac{h_{T_0}}{\rho_{T_0}} h_{T_0} | v |_{W^{1,p}(T_0)^d}.
\end{align*}
As geometric conditions to obtain global interpolation error estimates and to prove that this estimate converges to zero as $h  \to 0$, the \textit{shape-regularity condition} is widely used and well known. If the shape-regularity condition is violated, that is, the simplex becomes too flat as $h_{T_0} \to 0$, the quantity $\frac{h_{T_0}}{\rho_{T_0}} h_{T_0}$ may diverge.

To overcome the difficulty, we considered new strategy (\cite[Section 3]{IshKobTsu21a} and \cite[Section 2]{IshKobTsu21c}) to use anisotropic mesh partitions. Using an intermediate simplex $T \subset \mathbb{R}^d$ imposing Condition \ref{cond1} or \ref{cond2} described in Section \ref{Affinedef}, we construct two affine mappings $\Phi_T: \widehat{T} \to T$ and $\Phi_{T_0}: T \to T_0$. We first define the affine mapping $\Phi_T: \widehat{T} \to T$ as
\begin{align*}
\displaystyle
\Phi_T: \widehat{T} \ni \hat{x} \mapsto x := \Phi_T(\hat{x}) := {A}_T \hat{x} \in  T,
\end{align*}
where ${A}_T \in \mathbb{R}^{d \times d}$ is an invertible matrix. Details are given in Section \ref{sec221}. We next define an affine mapping $\Phi_{T_0}: T \to T_0$  as
\begin{align*}
\displaystyle
\Phi_{T_0}: T \ni x \mapsto x^{(0)} := \Phi_{T_0}(x) := {A}_{T_0} x + b_{T_0} \in T_0,
\end{align*}
where $b_{T_0} \in \mathbb{R}^d$ and ${A}_{T_0} \in \mathbb{R}^{d \times d}$ is an invertible matrix. Details are given in Section \ref{sec222}. We define an affine mapping $\Phi: \widehat{T} \to T_0$ as
\begin{align*}
\displaystyle
%\Phi := {\Phi}_{T_0} \circ {\Phi}_T: \widehat{T} \to T_0, \ x^{(0)} := \Phi (\hat{x}) =  ({\Phi}_{T_0} \circ {\Phi}_T)(\hat{x}) = {A} \hat{x} + b_{T_0}, 
\Phi := {\Phi}_{T_0} \circ {\Phi}_T: \widehat{T} \ni \hat{x} \mapsto x^{(0)} := \Phi (\hat{x}) =  ({\Phi}_{T_0} \circ {\Phi}_T)(\hat{x}) = {A} \hat{x} + b_{T_0} \in T_0, 
\end{align*}
where ${A} := {A}_{T_0} {A}_T \in \mathbb{R}^{d \times d}$. In our strategy, we use the affine mapping $\Phi$ instead of the mapping $\Phi_0$.
%and The two elements $\widehat{T}$ and $T_0$ are then affine equivalent. 

\subsection{Reference elements} \label{reference}
We define reference elements $\widehat{T} \subset \mathbb{R}^d$.

\subsubsection*{Two-dimensional case} \label{reference2d}
Let $\widehat{T} \subset \mathbb{R}^2$ be a reference triangle with vertices $\hat{x}_1 := (0,0)^T$, $\hat{x}_2 := (1,0)^T$, and $\hat{x}_3 := (0,1)^T$. 

\subsubsection*{Three-dimensional case} \label{reference3d}
In the three-dimensional case, we need to consider the following two cases (\roman{sone}) and (\roman{stwo}); also see Condition \ref{cond2}.

Let $\widehat{T}_1$ and $\widehat{T}_2$ be reference tetrahedra with the following vertices:
\begin{description}
   \item[(\roman{sone})] $\widehat{T}_1$ has the vertices $\hat{x}_1 := (0,0,0)^T$, $\hat{x}_2 := (1,0,0)^T$, $\hat{x}_3 := (0,1,0)^T$, and $\hat{x}_4 := (0,0,1)^T$;
 \item[(\roman{stwo})] $\widehat{T}_2$ has the vertices $\hat{x}_1 := (0,0,0)^T$, $\hat{x}_2 := (1,0,0)^T$, $\hat{x}_3 := (1,1,0)^T$, and $\hat{x}_4 := (0,0,1)^T$.
\end{description}
We thus set $\widehat{T} \in \{ \widehat{T}_1 , \widehat{T}_2 \}$.

\subsection{Affine mappings} \label{Affinedef}

\subsubsection{Construct of the mapping $\Phi_T: \widehat{T} \to T$} \label{sec221}
We define the affine mapping $\Phi_T: \widehat{T} \to T$ as
\begin{align}
\displaystyle
\Phi_T: \widehat{T} \ni \hat{x} \mapsto x := \Phi_T(\hat{x}) := {A}_T \hat{x} \in  T \label{mesh4}
\end{align}
with the element $T$ satisfying  Condition \ref{cond1} when $d=2$ or Condition \ref{cond2} when $d=3$, which are described later. We define the matrix $ {A}_T \in \mathbb{R}^{d \times d}$ as follows. We first define the diagonal matrix;
\begin{align}
\displaystyle
\widehat{A} :=  \diag (\alpha_1,\ldots,\alpha_d), \quad \alpha_i \in \mathbb{R}_+ \quad \forall i.  \label{mesh1}
\end{align}
When $d=2$, we define the regular matrix $\widetilde{A} \in \mathbb{R}^{2 \times 2}$ as
\begin{align}
\displaystyle
\widetilde{A} :=
\begin{pmatrix}
1 & s \\
0 & t \\
\end{pmatrix}, \label{mesh2}
\end{align}
with parameters
\begin{align*}
\displaystyle
s^2 + t^2 = 1, \quad t \> 0.
\end{align*}
%Note that for convenience, we use $\widetilde{A} := \widetilde{A}_1$. 
For the reference element $\widehat{T}$, let $\mathfrak{T}^{(2)}$ be the family of triangles
\begin{align*}
\displaystyle
T &= \Phi_T(\widehat{T}) = {A}_T (\widehat{T}), \quad {A}_{T} := \widetilde {A} \widehat{A}
\end{align*}
with vertices $x_1 := (0,0)^T$, $x_2 := (\alpha_1,0)^T$, and $x_3 :=(\alpha_2 s , \alpha_2 t)^T$. Then,  $\alpha_1 = |x_1 - x_2| \> 0$ and $\alpha_2 = |x_1 - x_3| \> 0$. 

When $d=3$, we define the regular matrices $\widetilde{A}_1, \widetilde{A}_2 \in \mathbb{R}^{3 \times 3}$ as
\begin{align}
\displaystyle
\widetilde{A}_1 :=
\begin{pmatrix}
1 & s_1 & s_{21} \\
0 & t_1  & s_{22}\\
0 & 0  & t_2\\
\end{pmatrix}, \
\widetilde{A}_2 :=
\begin{pmatrix}
1 & - s_1 & s_{21} \\
0 & t_1  & s_{22}\\
0 & 0  & t_2\\
\end{pmatrix} \label{mesh3}
\end{align}
with parameters
\begin{align*}
\displaystyle
\begin{cases}
s_1^2 + t_1^2 = 1, \ s_1 \> 0, \ t_1 \> 0, \ \alpha_2 s_1 \leq \alpha_1 / 2, \\
s_{21}^2 + s_{22}^2 + t_2^2 = 1, \ t_2 \> 0, \ \alpha_3 s_{21} \leq \alpha_1 / 2.
\end{cases}
\end{align*}
We thus set $\widetilde{A} \in \{ \widetilde{A}_1 , \widetilde{A}_2 \}$. For the reference elements $\widehat{T}_i$, $i=1,2$, let $\mathfrak{T}_i^{(3)}$, $i=1,2$, be the family of tetrahedra
\begin{align*}
\displaystyle
T_i &= \Phi_{T_i} (\widehat{T}_i) =  {A}_{T} (\widehat{T}_i), \quad {A}_{T} := \widetilde {A}_i \widehat{A}, \quad i=1,2,
\end{align*}
with vertices
\begin{align*}
\displaystyle
&x_1 := (0,0,0)^T, \ x_2 := (\alpha_1,0,0)^T, \ x_4 := (\alpha_3 s_{21}, \alpha_3 s_{22}, \alpha_3 t_2)^T, \\
&\begin{cases}
x_3 := (\alpha_2 s_1 , \alpha_2 t_1 , 0)^T \quad \text{for case (\roman{sone})}, \\
x_3 := (\alpha_1 - \alpha_2 s_1, \alpha_2 t_1,0)^T \quad \text{for case (\roman{stwo})}.
\end{cases}
\end{align*}
Then, $\alpha_1 = |x_1 - x_2| \> 0$, $\alpha_3 = |x_1 - x_4| \> 0$, and
\begin{align*}
\displaystyle
\alpha_2 =
\begin{cases}
|x_1 - x_3| \> 0  \quad \text{for case (\roman{sone})}, \\
|x_2 - x_3| \> 0  \quad \text{for case (\roman{stwo})}.
\end{cases}
\end{align*}

In the following, we impose conditions on $T \in \mathfrak{T}^{(2)}$ in the two-dimensional case and $T \in \mathfrak{T}_1^{(3)} \cup \mathfrak{T}_2^{(3)} =: \mathfrak{T}^{(3)}$ in the three-dimensional case. 

\begin{cond}[Case in which $d=2$] \label{cond1}
Let $T \in \mathfrak{T}^{(2)}$ with the vertices $x_i$ ($i=1,\ldots,3$) introduced in this section. We assume that $\overline{x_2 x_3}$ is the longest edge of $T$; i.e., $ h_T := |x_2 - x_ 3|$. Recall that $\alpha_1 = |x_1 - x_2|$ and $\alpha_2 = |x_1 - x_3|$. We then assume that $\alpha_2 \leq \alpha_1$. Note that $\alpha_1 = \mathcal{O}(h_T)$. 
\end{cond}
\begin{cond}[Case in which $d=3$] \label{cond2}
Let $T \in \mathfrak{T}^{(3)}$ with the vertices $x_i$ ($i=1,\ldots,4$) introduced in this section. Let $L_i$ ($1 \leq i \leq 6$) be the edges of $T$. We denote by $L_{\min}$  the edge of $T$ that has the minimum length; i.e., $|L_{\min}| = \min_{1 \leq i \leq 6} |L_i|$. We set $\alpha_2 := |L_{\min}|$ and assume that 
\begin{align*}
\displaystyle
&\text{the end points of $L_{\min}$ are either $\{ x_1 , x_3\}$ or $\{ x_2 , x_3\}$}.
\end{align*}
Among the four edges that share an end point with $L_{\min}$, we take the longest edge $L^{({\min})}_{\max}$. Let $x_1$ and $x_2$ be the end points of edge $L^{({\min})}_{\max}$. We thus have that
\begin{align*}
\displaystyle
\alpha_1 = |L^{(\min)}_{\max}| = |x_1 - x_2|.
\end{align*}
Consider cutting $\mathbb{R}^3$ with the plane that contains the midpoint of edge $L^{(\min)}_{\max}$ and is perpendicular to the vector $x_1 - x_2$. We then have two cases: 
\begin{description}
  \item[(Type \roman{sone})] $x_3$ and $x_4$  belong to the same half-space;
  \item[(Type \roman{stwo})] $x_3$ and $x_4$  belong to different half-spaces.
\end{description}
In each case, we set
\begin{description}
  \item[(Type \roman{sone})] $x_1$ and $x_3$ as the end points of $L_{\min}$, that is, $\alpha_2 =  |x_1 - x_3| $;
  \item[(Type \roman{stwo})] $x_2$ and $x_3$ as the end points of $L_{\min}$, that is, $\alpha_2 =  |x_2 - x_3| $.
\end{description}
Finally, we set $\alpha_3 = |x_1 - x_4|$. Note that we implicitly assume that $x_1$ and $x_4$ belong to the same half-space. In addition, note that $\alpha_1 = \mathcal{O}(h_T)$.
\end{cond}

\subsubsection{Construct of the mapping $\Phi_{T_0}: T \to T_0$}  \label{sec222}
Let $\Phi_{T_0}$ be the affine mapping defined as
\begin{align}
\displaystyle
\Phi_{T_0}: T \ni x \mapsto x^{(0)} := \Phi_{T_0}(x) := {A}_{T_0} x + b_{T_0} \in T_0, \label{affine10}
\end{align}
where $b_{T_0} \in \mathbb{R}^d$ and ${A}_{T_0} \in O(d)$ is a rotation and mirror imaging matrix. Note that none of the lengths of the edges of a simplex or the measure of the simplex is changed by the transformation.	

\begin{Ex}
As examples, we define the matrices $A_{T_0}$ as 
\begin{align*}
\displaystyle
A_{T_0} := 
\begin{pmatrix}
\cos \theta  & - \sin \theta \\
 \sin \theta & \cos \theta
\end{pmatrix}, \quad 
{A}_{T_0} := 
\begin{pmatrix}
 \cos \theta  & - \sin \theta & 0\\
 \sin \theta & \cos \theta & 0 \\
 0 & 0 & 1 \\
\end{pmatrix},
\end{align*}
where $\theta$ is an angle. %concering the positive the $x$-axis.
\end{Ex}

\subsection{Piola transformations} \label{Pioladef}
The Piola transformation $\Psi := {\Psi}_{T_0} \circ {\Psi}_T : \mathcal{C}(\widehat{T})^d \to \mathcal{C}({T}_0)^d$ is defined as
\begin{align*}
\displaystyle
\Psi :  \mathcal{C}(\widehat{T})^d  &\to  \mathcal{C}({T}_0)^d \\
\hat{v} &\mapsto v_0(x) :=  \Psi_{T_0} \circ {\Psi}_T(\hat{v})(x) = \frac{1}{\det(A)} A \hat{v}(\hat{x}), \quad A = A_{T_0} A_T
\end{align*}
with	two Piola transformations:
\begin{align*}
\displaystyle
{\Psi}_T:  \mathcal{C}(\widehat{T})^d &\to  \mathcal{C}({T})^d \nonumber\\
\hat{v} &\mapsto {v}({x}) := {\Psi}_T(\hat{v})({x}) := \frac{1}{\det ({A}_T)} {A}_T  \hat{v}(\hat{x}), \\
{\Psi}_{T_0}:  \mathcal{C}({T})^d &\to  \mathcal{C}({T}_0)^d \nonumber\\
{v} &\mapsto {v}_0({x}^{(0)}) := {\Psi}_{T_0}({v})({x}^{(0)}) := \frac{1}{\det ({A}_{T_0})} {A}_{T_0} {v}({x}).
\end{align*}

\subsection{Additional notation and assumption} \label{addinot}
For convenience, we introduce two definitions.

\begin{defi}
We define a parameter $\mathscr{H}_i$, $i=1,\ldots,d$, as
\begin{align*}
\displaystyle
\begin{cases}
\mathscr{H}_1 := \alpha_1, \quad \mathscr{H}_2 := \alpha_2 t \quad \text{if $d=2$}, \\
\mathscr{H}_1 := \alpha_1, \quad \mathscr{H}_2 := \alpha_2 t_1, \quad \mathscr{H}_3 := \alpha_3 t_2 \quad \text{if $d=3$}.
\end{cases}
\end{align*}
For a multi-index $\beta = (\beta_1,\ldots,\beta_d) \in \mathbb{N}_0^d$, we use the following notation:
\begin{align*}
\displaystyle
\mathscr{H}^{\beta} := \mathscr{H}_1^{\beta_1} \cdots \mathscr{H}_d^{\beta_d}, \quad \mathscr{H}^{- \beta} := \mathscr{H}_1^{- \beta_1} \cdots \mathscr{H}_d^{- \beta_d}.
\end{align*}
We also define  $\alpha^{\beta} :=  \alpha_{1}^{\beta_1} \cdots \alpha_{d}^{\beta_d}$ and $\alpha^{- \beta} :=  \alpha_{1}^{- \beta_1} \cdots \alpha_{d}^{- \beta_d}$. 
\end{defi}

\begin{figure}[tbhp]
\vspace{-8cm}
  \includegraphics[bb=0 0 1122 796,scale=0.55]{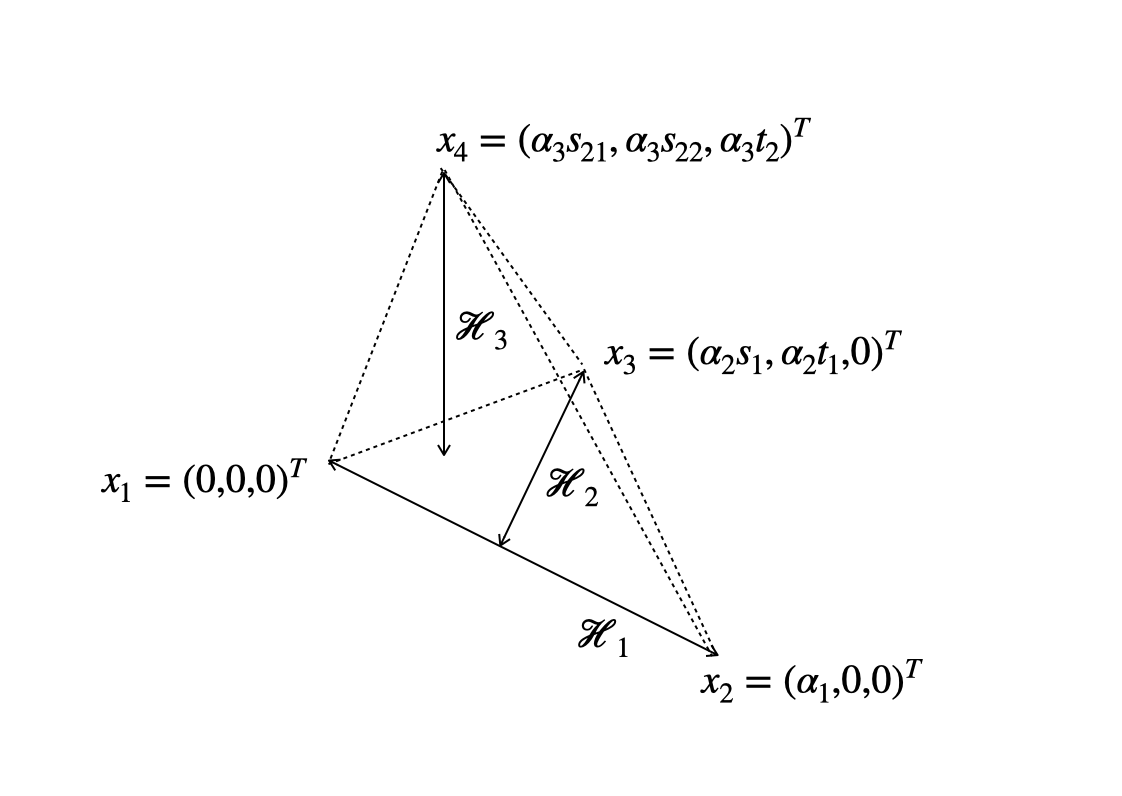}
\caption{New parameters $\mathscr{H}_i$, $i=1,2,3$}
\label{mathscrH}
\end{figure}

\begin{defi} \label{defi332}
We define vectors $r_n \in \mathbb{R}^d$, $n=1,\ldots,d$, as follows. If $d=2$,
\begin{align*}
\displaystyle
r_1 := (1 , 0)^T, \quad r_2 := (s,t)^T,
\end{align*}
and if $d=3$,
\begin{align*}
\displaystyle
&r_1 := (1 , 0,0)^T, \quad r_3 := ( s_{21}, s_{22} , t_2)^T, \\
&\begin{cases}
r_2 := ( s_1 ,  t_1 , 0)^T \quad \text{for case (\roman{sone})}, \\
r_2 := (- s_1,  t_1,0)^T \quad \text{for case (\roman{stwo})}.
\end{cases}
\end{align*}
For a sufficiently smooth function $\varphi_0$ and vector function $v_0 := (v_{0,1},\ldots,v_{0,d})^T$, we define the directional derivative as, for $i \in \{ 1 : d \}$,
\begin{align*}
\displaystyle
\frac{\partial \varphi_0}{\partial {r_i}^{(0)}} &:= [ ({A}_{T_0} r_i) \cdot  \nabla_{x^{(0)}} ] \varphi_0 = \sum_{i_0=1}^d ({A}_{T_0} r_i)_{i_0} \frac{\partial \varphi_0}{\partial x_{i_0}^{}} =  \sum_{i_0,j_0=1}^d [A_{T_0}]_{i_0 j_0} (r_i)_{j_0} \frac{\partial \varphi_0}{\partial x_{i_0}^{}}, \\
\frac{\partial v}{\partial r_i^{(0)}} &:= \left(\frac{\partial v_{0,1}}{\partial r_i^{(0)}}, \ldots, \frac{\partial v_{0,d}}{\partial r_i^{(0)}} \right)^T \\
&= ( [({A}_{T_0} r_i)  \cdot \nabla_{x^{(0)}}] v_{0,1}, \ldots, [ ({A}_{T_0} r_i)  \cdot \nabla_{x^{(0)}} ] v_{0,d} )^T,
\end{align*}
where ${A}_{T_0} \in O(d)$ is the orthogonal matrix defined in \eqref{affine10}. For a multi-index $\beta = (\beta_1,\ldots,\beta_d) \in \mathbb{N}_0^d$, we use the notation
\begin{align*}
\displaystyle
\partial^{\beta}_{r^{(0)}} \varphi_0 := \frac{\partial^{|\beta|} \varphi_0}{( \partial r_1^{(0)})^{\beta_1} \ldots ( \partial r_d^{(0)})^{\beta_d}}.
\end{align*}
Furthermore, for $\varphi = \varphi_0 \circ \Phi_{T_0}$ and $v = \Psi_{T_0}^{-1} v_0$, we define the directional derivative as, for $i \in \{ 1 : d \}$,
\begin{align*}
\displaystyle
\frac{\partial \varphi}{\partial {r_i}} &:= [ r_i \cdot  \nabla_{x} ] \varphi = \sum_{i_0=1}^d  (r_i)_{i_0} \frac{\partial \varphi}{\partial x_{i_0}}, \\
\frac{\partial v}{\partial r_i} &:= \left(\frac{\partial v_1}{\partial r_i}, \ldots, \frac{\partial v_d}{\partial r_i} \right)^T := ( [r_i  \cdot \nabla_{x}] v_1, \ldots, [  r_i  \cdot \nabla_{x} ] v_d )^T.
\end{align*}
For a multi-index $\beta = (\beta_1,\ldots,\beta_d) \in \mathbb{N}_0^d$, we use the notation
\begin{align*}
\displaystyle
\partial^{\beta}_{r} \varphi := \frac{\partial^{|\beta|} \varphi}{ \partial r_1^{\beta_1} \ldots \partial r_d^{\beta_d}}.
\end{align*}
\end{defi}

\begin{figure}[tbhp]
\vspace{-8cm}
  \includegraphics[bb=0 0 1290 796,scale=0.55]{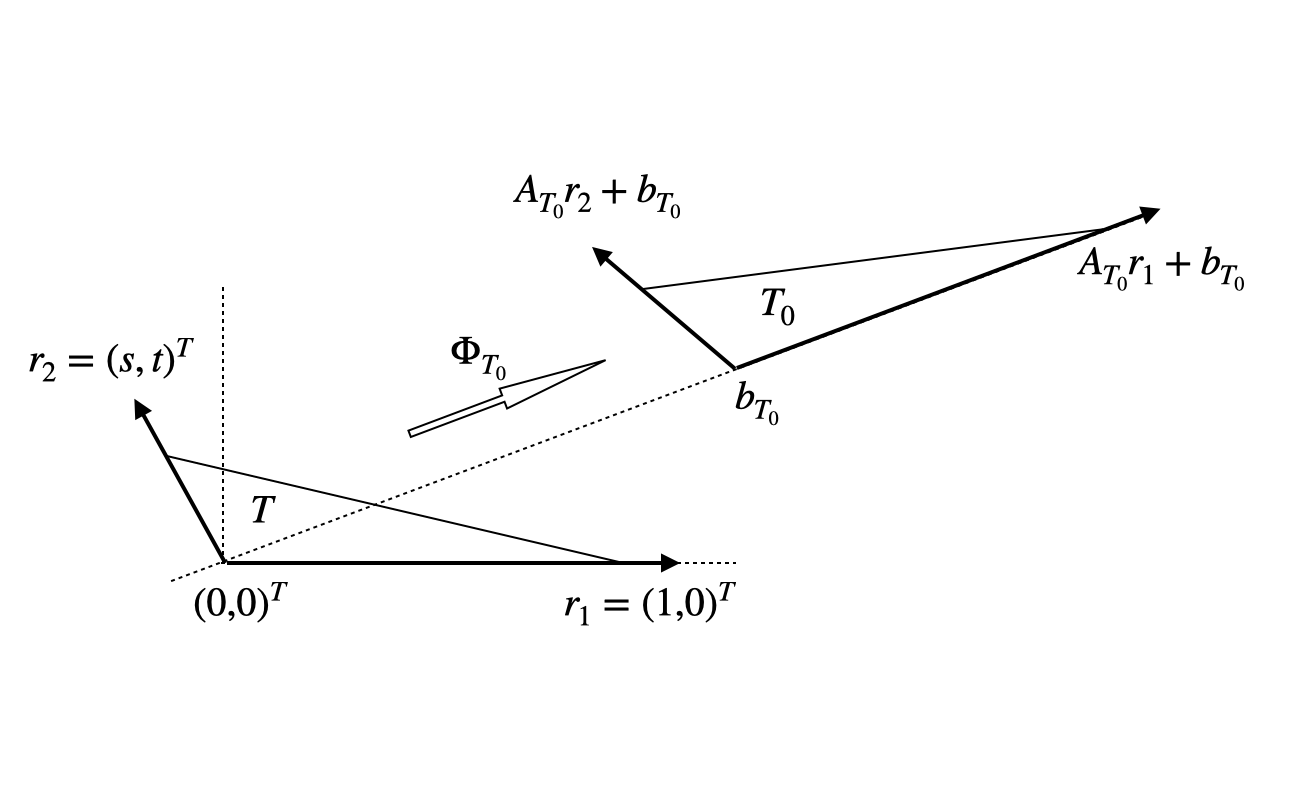}
\caption{Affine mapping $\Phi_{T_0}$ and vectors $r_i$, $i=1,2$}
\label{mathscrH}
\end{figure}

\begin{ass} \label{ass1}
In anisotropic interpolation error analysis, we may impose a geometric condition for the simplex $T$:
\begin{enumerate}
 \item If $d=2$, there are no additional conditions;
 \item If $d=3$, there exists a positive constant $M$ independent of $h_T$ such that $|s_{22}| \leq M \frac{\alpha_2 t_1}{\alpha_3}$. Note that if $s_{22} \neq 0$, this condition means that the order with respect to $h_T$ of $\alpha_3$ coincides with the order of $\alpha_2$, and if $s_{22} = 0$, the order of $\alpha_3$ may be different from that of $\alpha_2$. 
 %If $d=3$, we impose $s_{22} = 0$ or there exists a positive constant $m$ independent of $h_T$ such that $|s_{22}| \leq m \frac{h_2 t_1}{h_3}$.
%there exist positive constant $m_1$ and $m_2$ independent of $h_T$ such that $m_1 h_2 t_1 \leq h_3 s_{22}  \leq m_2 h_2 t_1$.
\end{enumerate}
\end{ass}

\begin{Note} \label{defi=theta}
Recall that
\begin{align*}
\displaystyle
&|s| \leq 1, \ \alpha_2 \leq \alpha_1 \quad \text{if $d=2$},\\
&|s_1|\leq 1, \ |s_{21}| \leq 1, \  \alpha_2 \leq \alpha_3\leq \alpha_1 \quad \text{if $d=3$}.
\end{align*}
When $d=3$, if Assumption \ref{ass1} is imposed, there exists a positive constant $M$ independent of $h_T$ such that $|s_{22}| \leq M \frac{\alpha_2 t_1}{\alpha_3}$. We thus have, if $d=2$,
\begin{align*}
\displaystyle
&\alpha_1 | [\widetilde{A}]_{j1} | \leq \mathscr{H}_j, \quad \alpha_2 | [\widetilde{A}]_{j2} | \leq \mathscr{H}_j, \quad j=1,2,
\end{align*}
and, if $d=3$, for $\widetilde{A} \in \{ \widetilde{A}_1 , \widetilde{A}_2  \}$,
\begin{align*}
\displaystyle
&\alpha_1 | [\widetilde{A}]_{j1} | \leq \mathscr{H}_j, \quad \alpha_2 | [\widetilde{A}]_{j2} | \leq \mathscr{H}_j, \quad \alpha_3 | [\widetilde{A}]_{j3} | \leq \max \{ 1,M\} \mathscr{H}_j,  \quad j=1,2,3.
\end{align*}
\end{Note}

\subsection{New parameters}
 We proposed two geometric parameters, $H_{T}$ and $H_{T_0}$, in \cite{IshKobTsu21a}.
 \begin{defi} \label{defi1}
 The parameter $H_T$ is defined as
\begin{align*}
\displaystyle
H_T := \frac{\prod_{i=1}^d \alpha_i}{|T|} h_T,
\end{align*}
and the parameter $H_{T_0}$ is defined as
\begin{align*}
\displaystyle
H_{T_0} := \frac{h_{T_0}^2}{|T_0|} \min_{1 \leq i \leq 3} |L_i|  \quad \text{if $d=2$}, \quad H_{T_0} := \frac{h_{T_0}^2}{|T_0|} \min_{1 \leq i , j \leq 6, i \neq j} |L_i| |L_j| \quad \text{if $d=3$}
\end{align*}
where $L_i$ denotes the edges of the simplex $T_0 \subset \mathbb{R}^d$. 
\end{defi}

The following lemma shows the equivalence between $H_T$ and $H_{T_0}$.

 \begin{lem} \label{lem2}
 It holds that
\begin{align*}
\displaystyle
\frac{1}{2} H_{T_0} \< H_T  \< 2 H_{T_0}.
\end{align*}
Furthermore, in the two-dimensional case, $H_{T_0}$ is equivalent to the circumradius $R_2$ of $T_0$.
 \end{lem}
 
 \begin{pf*}
The proof is found in \cite[Lemma 3]{IshKobTsu21a}.
\qed
\end{pf*}

We introduce the geometric condition proposed in \cite{IshKobTsu21a}, which is equivalent to the maximum-angle condition \cite{IshKobSuzTsu21d}.

\begin{ass} \label{ass3}
A family of meshes $\{ \mathbb{T}_h\}$ has a semi-regular property if there exists $\gamma_0 \> 0$ such that \eqref{NewGeo}. 
%\begin{align}
%\displaystyle
%\frac{H_{T_0}}{h_{T_0}} \leq \gamma_0 \quad \forall \mathbb{T}_h \in \{ \mathbb{T}_h \}, \quad \forall T_0 \in \mathbb{T}_h. \label{NewGeo}
%\end{align}
Equivalently, there exists $\gamma_1 \> 0$ such that
\begin{align}
\displaystyle
\frac{H_{T}}{h_{T}} \leq \gamma_1 \quad \forall \mathbb{T}_h \in \{ \mathbb{T}_h \}, \quad \forall T_0 \in \mathbb{T}_h, \quad T = \Phi_{T_0}^{-1} (T_0). \label{NewGeo2}
\end{align}
\end{ass}

\begin{Rem} \label{goodel}
In \cite{IshKobTsu21c}, we considered good elements on meshes. %"Good elements" on meshes mean that there exists a positive constant $\gamma_0 \> 0$ satisfying \eqref{NewGeo}. 
On anisotropic meshes, the good elements may satisfy conditions such as
\begin{description}
  \item[($d=2$)] $\alpha_2 \approx \alpha_2 t = \mathscr{H}_2$;
  \item[($d=3$)] $\alpha_2 \approx \alpha_2 t_1 = \mathscr{H}_2$ and $\alpha_3 \approx \alpha_3 t_2 = \mathscr{H}_3$.
\end{description}
\end{Rem}

We have the following theorem concerning the new condition.

\begin{thr} \label{thr3}
Condition \eqref{NewGeo} holds if and only if there exist $0 \< \gamma_2, \gamma_3 \< \pi$ such that 
\begin{align}
\displaystyle
d=2: \quad \theta_{T_0,\max} \leq \gamma_2 \quad \forall \mathbb{T}_h \in \{ \mathbb{T}_h \}, \quad \forall T_0 \in \mathbb{T}_h, \label{M.A.C.2}
\end{align}
where $\theta_{T_0,\max}$ is the maximum angle of $T_0$, and 
\begin{align}
\displaystyle
d=3: \quad \theta_{T_0, \max} \leq \gamma_3, \quad \psi_{T_0, \max} \leq \gamma_3 \quad \forall \mathbb{T}_h \in \{ \mathbb{T}_h \}, \quad \forall T_0 \in \mathbb{T}_h,  \label{M.A.C.3}
\end{align}
where $\theta_{T_0, \max}$ is the maximum angle of all triangular faces of the tetrahedron $T_0$ and $\psi_{T_0, \max}$ is the maximum dihedral angle of $T_0$. Conditions \eqref{M.A.C.2} and \eqref{M.A.C.3} are called the \textit{maximum-angle condition}.
\end{thr}

\begin{pf*}
In the case that $d=2$, we use the previous result presented in \cite{Kri91}; i.e., there exists a constant $\gamma_3  \> 0$ such that
\begin{align*}
\displaystyle
\frac{R_2}{h_{T_0}} \leq \gamma_3 \quad \forall \mathbb{T}_h \in \{ \mathbb{T}_h \}, \quad \forall T_0 \in \mathbb{T}_h,
\end{align*}
if and only if condition \eqref{M.A.C.2} is satisfied. Combining this result with the fact that $H_{T_0}$ is equivalent to the circumradius $R_2$ of $T_0$, we have the desired conclusion. In the case that $d=3$, the proof can be found in the recent paper \cite{IshKobSuzTsu21d}.
\qed	
\end{pf*}

\begin{lem}[Euclidean condition number] \label{lem351}
It holds that
\begin{subequations} \label{CN331}
\begin{align}
\displaystyle
\| \widehat{{A}} \|_2 &\leq  h_{T}, \quad \| \widehat{{A}} \|_2 \| \widehat{{A}}^{-1} \|_2 = \frac{\max \{\alpha_1 , \cdots, \alpha_d \}}{\min \{\alpha_1 , \cdots, \alpha_d \}}, \label{CN331a} \\
\| \widetilde{{A}} \|_2 &\leq 
\begin{cases}
\sqrt{2} \quad \text{if $d=2$}, \\
2  \quad \text{if $d=3$},
\end{cases}
\quad \| \widetilde{{A}} \|_2 \| \widetilde{{A}}^{-1} \|_2 \leq
\begin{cases}
\frac{\alpha_1 \alpha_2}{|T|} = \frac{H_{T}}{h_{T}} \quad \text{if $d=2$}, \\
\frac{2}{3} \frac{\alpha_1 \alpha_2 \alpha_3}{|T|} = \frac{2}{3} \frac{H_{T}}{h_{T}} \quad \text{if $d=3$},
\end{cases} \label{CN331b} \\
\| A_{T_0} \|_2 &= 1, \quad \| A_{T_0}^{-1} \|_2 = 1. \label{CN331c}
\end{align}
\end{subequations}
Furthermore, we have
\begin{align}
\displaystyle
| \det ({A}_T) | = | \det(\widetilde{{A}}) | | \det (\widehat{{A}}) | = d ! |T |, \quad | \det ({A}_{T_0}) | = 1. \label{CN332}
\end{align}
\end{lem}

\begin{pf*}
The proof of  \eqref{CN331b} is found in \cite[(4.4), (4.5), (4.6) and (4.7)]{IshKobTsu21a}. The inequality \eqref{CN331a} is easily proved.  Because $A_{T_0} \in O(d)$, one can easily have $A_{T_0}^{-1} \in O(d)$ and \eqref{CN331c}. The proof of equality \eqref{CN332} is standard.
\qed
\end{pf*}

\subsection{Raviart--Thomas finite element generation} \label{FEGeneration}
We follow the procedure described in \cite[Section 3.5]{IshKobTsu21a} and \cite[Section 1.4.1 and 1.2.1]{ErnGue04}.

%\subsection{Raviart--Thomas finite element method} \label{RTSpace}
For a simplex $T_0 \subset \mathbb{R}^d$, we define the local Raviart--Thomas polynomial space of order $k \in \mathbb{N}_0$ as follows:
\begin{align}
\displaystyle
RT^k(T_0 ) := \mathcal{P}^k(T_0 )^d + x \mathcal{P}^k(T_0 ), \quad x \in \mathbb{R}^d. \label{RTp}
\end{align}

Let $\widehat{T} \subset \mathbb{R}^d$ be the reference element defined in Sections \ref{reference}.  The Raviart--Thomas finite element on the reference element is defined by the triple $\{ \widehat{T} , \widehat{P} , \widehat{\Sigma\}}$ as follows:
 \begin{enumerate}
 \item $\widehat{P} := RT^k(\widehat{T})$;
 \item setting $N^{(RT)} := \dim RT^k$, $\widehat{\Sigma}$ is a set $\{ \hat{\chi}_{i} \}_{1 \leq i \leq N^{(RT)}}$ of $N^{(RT)}$ linear forms with its components such that, for any $\hat{q} \in \widehat{P}$, 
\begin{align}
\displaystyle
%\chi^1_{i}(v) &:= \int_{F_{i}} v \cdot n_{F_i} p_j ds, \quad \forall p_j \in \mathcal{P}^k(F_i), \quad F_i \subset \partial T, \label{pre5}
& \int_{\widehat{F}} \hat{q} \cdot \hat{n}_{\widehat{F}} \hat{p}_k d\hat{s}, \quad \forall \hat{p}_k \in \mathcal{P}^k(\widehat{F}), \quad \widehat{F} \subset \partial \widehat{T}, \label{RTdof1} \\
 &\int_{\widehat{T}} \hat{q} \cdot \hat{p}_{k-1} d\hat{x}, \quad \forall \hat{p}_{k-1} \in \mathcal{P}^{k-1}(\widehat{T})^d, \label{RTdof2}
\end{align}
where $\hat{n}_{\widehat{F}}$ denotes the outer unit normal vector of $\widehat{T}$ on the face $\widehat{F}$. Note that for $k=0$, the local degrees of freedom of type \eqref{RTdof2} are violated.
 \end{enumerate}
For the simplicial Raviart--Thomas element in $\mathbb{R}^d$, it holds that
\begin{align}
\displaystyle
\dim RT^k(\widehat{T}) =
\begin{cases}
(k+1)(k+3) \quad \text{if $d=2$}, \\
\frac{1}{2}(k+1)(k+2)(k+4) \quad \text{if $d=3$}.
\end{cases} \label{RTdof}
\end{align}
%It is known that the Raviart--Thomas finite element with the set of linear forms $\Sigma := \{ \chi^{1}_{i,j}, \chi^2_{\ell} \}$ is unisolvent; e.g., see \cite[Proposition 2.3.4]{BofBreFor13}. The triple $\{ T , RT^k , \Sigma \}$ is then a finite element.
The Raviart--Thomas finite element with the local degrees of freedom of \eqref{RTdof1} and \eqref{RTdof2} is unisolvent; for example, see \cite[Proposition 2.3.4]{BofBreFor13}. 
%The triple $\{ \widehat{T} , \widehat{P} , \widehat{\Sigma} \}$ is then a finite element.
 
 We set the domain of the local Raviart--Thomas interpolation to $V(\widehat{T}) := W^{1,1}(\widehat{T})^d$ ; for example, see \cite[p. 188]{ErnGue21a}. The local Raviart--Thomas interpolation $I_{\widehat{T}}^{RT^k}: V(\widehat{T}) \to \widehat{P}$ is then defined as follows: For any $\hat{v} \in V(\widehat{T})$,
 \begin{align}
\displaystyle
\int_{\widehat{F}} ( I_{\widehat{T}}^{RT^k} \hat{v} - \hat{v} ) \cdot \hat{n}_{\widehat{F}} \hat{p}_k d\hat{s} &= 0 \quad \forall \hat{p}_k \in \mathcal{P}^k(\widehat{F}), \quad \widehat{F} \subset \partial \widehat{T}, \label{RTinter1}
\end{align}
and if $k \geq 1$,
\begin{align}
\displaystyle
\int_{\widehat{T}} ( I_{\widehat{T}}^{RT^k} \hat{v} - \hat{v} ) \cdot \hat{q}_{k-1} d\hat{x} &= 0 \quad \forall \hat{q}_{k-1} \in \mathcal{P}^{k-1}(\widehat{T})^d. \label{RTinter2}
\end{align}

The triples $\{ {T}, {P}, {\Sigma}\}$ and $\{ {T}_0, {P}_0, {\Sigma}_0 \}$ are defined as
\begin{align*}
\displaystyle
\begin{cases}
\displaystyle
{T} = {\Phi_T}(\widehat{T}); \\
\displaystyle
P = \{ {\Psi_T}(\hat{p}) ; \ \hat{p} \in \widehat{P} \}; \\
\displaystyle
{\Sigma} = \{ \{ {\chi}_{i} \}_{1 \leq i \leq N^{(RT)}}; \ {\chi}_{i} = \hat{\chi}_i({\Psi_T}^{-1}({p})), \forall {p} \in P, \  \hat{\chi}_i \in \widehat{\Sigma}  \},
\end{cases}
\end{align*}
and
\begin{align*}
\displaystyle
\begin{cases}
\displaystyle
T_0 = {\Phi}_{T_0}({T}); \\
\displaystyle
P_0 = \{ {\Psi}_{T_0}({p}) ; \ {p} \in P) \}; \\
\displaystyle
{\Sigma}_{0} = \{ \{ {\chi}_{T_0,i} \}_{1 \leq i \leq N^{(RT)}}; \ {\chi}_{T_0,i} = {\chi}_i( {\Psi}_{T_0}^{-1}({p}_0)), \forall {p}_0 \in P_0, {\chi}_i \in {\Sigma}  \}.
\end{cases}
\end{align*}
The triple $\{ {T}, {P}, {\Sigma}\}$ and $\{ {T}_0, {P}_0, {\Sigma}_0 \}$ are then the Raviart--Thomas finite elements. Furthermore, let 
\begin{align}
\displaystyle
I_{{T}}^{RT^k}: V(T) \to P, \quad I_{{T}_0}^{RT^k}: V(T_0) \to P_0 \label{RTlocal}
\end{align}
be the associated local Raviart--Thomas interpolations defined in \eqref{RTinter1} and \eqref{RTinter2}.

\section{Scaling argument} \label{scaling=sec}
We present estimates related to the scaling argument corresponding to \cite[Lemma 1.113]{ErnGue04}. The estimates play significant roles in our analysis. %Calculated results which are used in proofs are found in \cite{phd}.

\begin{lem}
Let $p \in [1,\infty)$. Let $T \in \mathfrak{T}^{(d)}$ satisfy Condition \ref{cond1} or Condition \ref{cond2}. Let $T_0 \subset \mathbb{R}^d$ be a simplex such that $T = \Phi_{T_0}^{-1}(T_0)$. It holds that, for any $\hat{v} = (\hat{v}_1,\ldots,\hat{v}_d)^T \in L^{p}(\widehat{T})^d$ with ${v} = ({v}_1,\ldots,{v}_d)^T := {\Psi}_{T} \hat{v}$ and ${v}_0 = ({v}_{0,1},\ldots,{v}_{0,d})^T := {\Psi}_{T_0} {v}$, 
\begin{align}
\displaystyle
\| {v}_0 \|_{L^p({T_0})^d} \leq c |\det({A}_T)|^{\frac{1-p }{p}} \| \widetilde{{A}} \|_{2} \left(  \sum_{j=1}^d \alpha_j^p \| \hat{v}_j \|_{L^p(\widehat{T})}^p \right)^{1/p}. \label{RT41}
\end{align}

Let $\ell , m \in \mathbb{N}_0$ and  $k \in \mathbb{N}$ with $1 \leq k \leq d$. Let $\beta := (\beta_1,\ldots,\beta_d) \in \mathbb{N}_0^d$ and $\gamma := (\gamma_1,\ldots,\gamma_d) \in \mathbb{N}_0^d$ be multi-indices with $|\beta| = \ell$ and $|\gamma| = m$, respectively. It then holds that, for any $\hat{v} = (\hat{v}_1,\ldots,\hat{v}_d)^T \in W^{ |\beta|+|\gamma|,p}(\widehat{T})^d$ with ${v} = ({v}_1,\ldots, {v}_d)^T := {\Psi}_{T} \hat{v}$ and ${v}_0 = ({v}_{0,1},\ldots,{v}_{0,d})^T := {\Psi}_{T_0} {v}$, 
\begin{align}
\displaystyle
\left \|  \partial_{\hat{x}}^{\beta}  \partial_{\hat{x}}^{\gamma} \hat{v}_k \right \|_{L^{p}(\widehat{T})} &\leq c  |\det ({A}_T)|^{\frac{p-1}{p}} \alpha_k^{-1} \| \widetilde{{A}}^{-1} \|_{2}  \sum_{|\varepsilon| = |\beta|+|\gamma|} \alpha^{\varepsilon} \left \| \partial_{r^{(0)}}^{ \varepsilon} v_0  \right \|_{L^p(T_0)^d}. \label{RT42}
\end{align}
If Assumption \ref{ass1} is imposed, it holds that, for any $\hat{v} = (\hat{v}_1,\ldots,\hat{v}_d)^T \in W^{|\beta|+|\gamma|,p}(\widehat{T})^d$ with ${v} = ({v}_1,\ldots,{v}_d)^T := {\Psi}_T \hat{v}$ and ${v}_0 = (v_{0,1}, \ldots,v_{0,d})^T := \Psi_{T_0}{v}$, 
\begin{align}
\displaystyle
\left \|  \partial_{\hat{x}}^{\beta} \partial_{\hat{x}}^{\gamma} \hat{v}_k \right \|_{L^{p}(\widehat{T})} \leq c  |\det ({A}_T)|^{\frac{p-1}{p}} \alpha_k^{-1} \| \widetilde{{A}}^{-1} \|_{2} \sum_{|\varepsilon| = |\beta|+|\gamma|} \mathscr{H}^{\varepsilon} \| \partial_x^{\varepsilon} (\Psi_{T_0}^{-1}  v_0 ) \|_{L^p( \Phi_{T_0}^{-1} (T_0))^d}.  \label{RT43}
\end{align}
%Here, for $p = \infty$ and any positive real $x$, $x^{- \frac{1}{p}} = 1$.
\end{lem}

\begin{pf*}
Because the space $\mathcal{C}(\widehat{T})^d$ is dense in the space $L^{p}(\widehat{T})^d$, we show \eqref{RT41} for $\hat{v} \in \mathcal{C}(\widehat{T})^d$ with ${v} = {\Psi}_T \hat{v}$ and  ${v}_0 = {\Psi}_{T_0} {v}$. The following inequality is found in \cite[Lemma 1.113]{ErnGue04}. There exists a positive constant $c$ such that
\begin{align*}
\displaystyle
\| v_0 \|_{L^{p}(T_0)^d} &\leq c \| {A}_{T_0} \|_2 |\det ({A}_{T_0})|^{ -\frac{1}{p^{\prime}}} \| v \|_{L^{p}(T)^d}, 
\end{align*}
which leads to, using \eqref{CN331c} and \eqref{CN332},
\begin{align}
\displaystyle
\| v_0 \|_{L^{p}(T_0)^d} &\leq c \| v \|_{L^{p}(T)^d}, \label{newRT56}
\end{align}
where $p^{\prime}$ is a real number such that $\frac{1}{p} + \frac{1}{p^{\prime}} = 1$. From the definition of the Piola transformation, for $i=1,\ldots,d$,
\begin{align*}
\displaystyle
{v}_i({x})  &= \frac{1}{\det ({A}_T)} \sum_{j=1}^d \widetilde{{A}}_{ij} \alpha_j \hat{v}_j (\hat{x}).
\end{align*}
If $1 \leq p \< \infty$, for $i=1,\ldots,d$,
\begin{align*}
\displaystyle
\| v \|^p_{L^p(T)^d} &= \sum_{i=1}^d \| v_i \|_{L^p(T)}^p \leq c |\det({A}_T)|^{1-p} \| \widetilde{{A}} \|_2^p \sum_{j=1}^d  \alpha_{j}^p \| \hat{v}_j \|^p_{L^p(\widehat{T})},
\end{align*}
which leads to \eqref{RT41} together with \eqref{newRT56}. %\eqref{Anorm} and \eqref{CN332}. 

Because the space $\mathcal{C}^{\ell+m}(\widehat{T})^d$ is dense
in the space $W^{\ell+m,p}(\widehat{T})^d$,  we prove \eqref{RT42} and \eqref{RT43} for $\hat{v} \in \mathcal{C}^{\ell+m}(\widehat{T})^d$ with ${v} = {\Psi}_T \hat{v}$ and  ${v}_0 = {\Psi}_{T_0} {v}$. Using \eqref{Anorm}, \eqref{CN331c} and \eqref{CN332}, through a simple calculation, we have, for $1 \leq k \leq d$, 
\begin{align*}
\displaystyle
| \partial^{\beta + \gamma}_{\hat{x}}  \hat{v}_k| &= \left| \frac{\partial^{|\beta| + |\gamma|}}{\partial \hat{x}_1^{\beta_1} \cdots \partial \hat{x}_d^{\beta_d} \partial \hat{x}_1^{\gamma_1} \cdots \partial \hat{x}_d^{\gamma_d}}  \hat{v}_k \right| \notag\\
&\hspace{-1.5cm}  \leq c | \det ({A}_T) | \| \widetilde{{A}}^{-1} \|_2 \alpha_k^{-1} \sum_{\nu=1}^d \alpha^{\beta} \alpha^{\gamma} \notag\\
&\hspace{-1.2cm} \Biggl | \underbrace{\sum_{i_1^{(1)},i_1^{(0,1)} = 1}^d  [{A}_{T_0}]_{i_1^{(0,1)} i_1^{(1)}} (r_1)_{i_1^{(1)}} \cdots   \sum_{i_{\beta_1}^{(1)},i_{\beta_1}^{(0,1)} = 1}^d  [{A}_{T_0}]_{i_{\beta_1}^{(0,1)} i_{\beta_1}^{(1)}} (r_1)_{i_{\beta_1}^{(1)}}  }_{\beta_1 \text{times}} \cdots \notag \\
&\hspace{-1.2cm} \underbrace{ \sum_{i_1^{(d)} , i_1^{(0,d)} = 1}^d   [{A}_{T_0}]_{i_1^{(0,d)} i_1^{(d)}} (r_d)_{i_{1}^{(d)}}  \cdots \sum_{i_{\beta_d}^{(d)} , i_{\beta_d}^{(0,d)}= 1}^d [{A}_{T_0}]_{i_{\beta_d}^{(0,d)} i_{\beta_d}^{(d)}} (r_d)_{i_{\beta_d}^{(d)}}  }_{\beta_d \text{times}}  \notag \\
&\hspace{-1.2cm}  \underbrace{\sum_{j_1^{(1)} , j_1^{(0,1)} = 1}^d   [{A}_{T_0}]_{j_1^{(0,1)} j_1^{(1)}} (r_1)_{j_{1}^{(1)}}   \cdots  \sum_{j_{\gamma_1}^{(1)} , j_{\gamma_1}^{(0,1)}= 1}^d  [{A}_{T_0}]_{j_{\gamma_1}^{(0,1)} j_{\gamma_1}^{(1)}} (r_1)_{ j_{\gamma_1}^{(1)}} }_{\gamma_1 \text{times}}  \cdots   \notag \\
&\hspace{-1.2cm}  \underbrace{\sum_{j_1^{(d)} , j_1^{(0,d)} = 1}^d  [{A}_{T_0}]_{j_1^{(0,d)} j_1^{(d)}} (r_d)_{j_1^{(d)}}  \cdots \sum_{j_{\gamma_d}^{(d)} , j_{\gamma_d}^{(0,d)}= 1}^d [{A}_{T_0}]_{j_{\gamma_d}^{(0,d)} j_{\gamma_d}^{(d)}} (r_d)_{j_{\gamma_d}^{(d)}} }_{\gamma_d \text{times}} \notag \\
&\hspace{-1.2cm}  \underbrace{\frac{\partial^{\beta_1}}{\partial {x}_{i_1^{(0,1)}}^{(0)} \cdots \partial {x}_{i_{\beta_1}^{(0,1)}}^{(0)}}}_{\beta_1 \text{times}} \cdots \underbrace{\frac{\partial^{\beta_d}}{\partial {x}_{i_1^{(0,d)}}^{(0)} \cdots \partial {x}_{i_{\beta_d}^{(0,d)}}^{(0)}}}_{\beta_d \text{times}}  \notag \\
&\hspace{-1.2cm} \underbrace{ \frac{\partial^{\gamma_1}}{ \partial {x}_{j_1^{(0,1)}}^{(0)} \cdots \partial {x}_{j_{\gamma_1}^{(0,1)}}^{(0)} } }_{\gamma_1 \text{times}} \cdots \underbrace{\frac{\partial^{\gamma_d}}{ \partial {x}_{j_1^{(0,d)}}^{(0)} \cdots \partial {x}_{j_{\gamma_d}^{(0,d)}}^{(0)} }}_{\gamma_d \text{times}} v_{0,\nu} \Biggr | \\
&\hspace{-1.5cm}  \leq c | \det ({A}_T) | \| \widetilde{{A}}^{-1} \|_2 \alpha_k^{-1} \sum_{\nu=1}^d \sum_{|\varepsilon| = |\beta|+|\gamma|} \alpha^{\varepsilon} | \partial_{r}^{\varepsilon} v_{\nu} |.
\end{align*}
Because $1 \leq p \< \infty$, it holds that, for $1 \leq k \leq d$,
\begin{align*}
\displaystyle
\left \|  \partial^{\beta}_{\hat{x}} \partial^{\gamma}_{\hat{x}}  \hat{v}_k  \right \|_{L^{p}(\widehat{T})}^p 
&\leq c  | \det ({A}_T) |^{p-1} \| \widetilde{{A}}^{-1} \|_2^p \alpha_k^{-p} \sum_{|\varepsilon| = |\beta|+|\gamma|} {\alpha}^{\varepsilon p} \int_{T_0} |\partial^{\varepsilon}_{r^{(0)}} v_0 |^p dx^{(0)},
\end{align*}
which leads to \eqref{RT42} together with \eqref{jensen}.

Using \eqref{Anorm} and Note \ref{defi=theta}, through a simple calculation, we have, for $1 \leq k \leq d$,
\begin{align*}
\displaystyle
| \partial^{\beta + \gamma}_{\hat{x}} \hat{v}_k |
&= \left|  \frac{\partial^{|\beta| + |\gamma|}}{\partial \hat{x}_1^{\beta_1} \cdots \partial \hat{x}_d^{\beta_d} \partial \hat{x}_1^{\gamma_1} \cdots \partial \hat{x}_d^{\gamma_d}}  \hat{v}_k \right| \\
&\hspace{-1.5cm}  \leq c | \det ({A}_T) |  \alpha_k^{-1} \sum_{\eta=1}^d | [\widetilde{{A}}^{-1}]_{k \eta}  | \\
&\hspace{-1.2cm}  \underbrace{\sum_{i_1^{(1)} = 1}^d \alpha_1 | [\widetilde{{A}}]_{i_1^{(1)} 1} | \cdots   \sum_{i_{\beta_1}^{(1)} = 1}^d \alpha_1 | [\widetilde{{A}}]_{i_{\beta_1}^{(1)} 1} |   }_{\beta_1 \text{times}} \cdots \underbrace{ \sum_{i_1^{(d)}  = 1}^d \alpha_d | [\widetilde{{A}}]_{i_1^{(d)} d} |  \cdots \sum_{i_{\beta_d}^{(d)} = 1}^d \alpha_d | [\widetilde{{A}}]_{i_{\beta_d}^{(d)} d} | }_{\beta_d \text{times}}  \notag \\
&\hspace{-1.2cm}  \underbrace{\sum_{j_1^{(1)}  = 1}^d \alpha_1 | [\widetilde{{A}}]_{j_1^{(1)} 1}|  \cdots  \sum_{j_{\gamma_1}^{(1)}= 1}^d \alpha_1 | [\widetilde{{A}}]_{j_{\gamma_1}^{(1)} 1} | }_{\gamma_1 \text{times}}  \cdots  \underbrace{\sum_{j_1^{(d)} = 1}^d \alpha_d | [\widetilde{{A}}]_{j_1^{(d)} d} | \cdots \sum_{j_{\gamma_d}^{(d)} = 1}^d \alpha_d | [\widetilde{{A}}]_{j_{\gamma_d}^{(d)} d} | }_{\gamma_d \text{times}} \notag \\
&\hspace{-1.2cm} \Biggl | \underbrace{\frac{\partial^{\beta_1}}{\partial {x}_{i_1^{(1)}} \cdots \partial {x}_{i_{\beta_1}^{(1)}}}}_{\beta_1 \text{times}} \cdots \underbrace{\frac{\partial^{\beta_d}}{\partial {x}_{i_1^{(d)}} \cdots \partial {x}_{i_{\beta_d}^{(d)}}}}_{\beta_d \text{times}} \underbrace{ \frac{\partial^{\gamma_1}}{ \partial {x}_{j_1^{(1)}} \cdots \partial {x}_{j_{\gamma_1}^{(1)}} } }_{\gamma_1 \text{times}} \cdots \underbrace{\frac{\partial^{\gamma_d}}{ \partial {x}_{j_1^{(d)}} \cdots \partial {x}_{j_{\gamma_d}^{(d)}} }}_{\gamma_d \text{times}} v^{}_{\eta} \Biggr|  \\
&\hspace{-1.5cm}  \leq c | \det ({A}_T) | \| \widetilde{{A}}^{-1} \|_2 \alpha_k^{-1} \\
&\hspace{-1.2cm} \sum_{\eta=1}^d   \underbrace{\sum_{i_1^{(1)} = 1}^d  \cdots   \sum_{i_{\beta_1}^{(1)} = 1}^d  }_{\beta_1 \text{times}} \cdots \underbrace{ \sum_{i_1^{(d)}  = 1}^d  \cdots \sum_{i_{\beta_d}^{(d)} = 1}^d  }_{\beta_d \text{times}}   \underbrace{\sum_{j_1^{(1)}  = 1}^d  \cdots  \sum_{j_{\gamma_1}^{(1)}= 1}^d  }_{\gamma_1 \text{times}}  \cdots  \underbrace{\sum_{j_1^{(d)} = 1}^d   \cdots \sum_{j_{\gamma_d}^{(d)} = 1}^d  }_{\gamma_d \text{times}} \notag \\
&\hspace{-1.2cm}  \underbrace{\mathscr{H}_{i_1^{(1)}} \cdots \mathscr{H}_{i_{\varepsilon_1}^{(1)}}}_{\beta_1 \text{times}} \cdots  \underbrace{ \mathscr{H}_{i_1^{(d)}} \cdots \mathscr{H}_{i_{\varepsilon_d}^{(d)}}}_{\beta_d \text{times}}  \underbrace{\mathscr{H}_{j_1^{(1)}} \cdots \mathscr{H}_{j_{\varepsilon_1}^{(1)}}}_{\gamma_1 \text{times}} \cdots  \underbrace{ \mathscr{H}_{j_1^{(d)}} \cdots \mathscr{H}_{j_{\varepsilon_d}^{(d)}}}_{\gamma_d \text{times}} \notag \\
&\hspace{-1.2cm} \Biggl |  \underbrace{\frac{\partial^{\beta_1}}{\partial {x}_{i_1^{(1)}} \cdots \partial {x}_{i_{\beta_1}^{(1)}}}}_{\beta_1 \text{times}} \cdots \underbrace{\frac{\partial^{\beta_d}}{\partial {x}_{i_1^{(d)}} \cdots \partial {x}_{i_{\beta_d}^{(d)}}}}_{\beta_d \text{times}} \underbrace{ \frac{\partial^{\gamma_1}}{ \partial {x}_{j_1^{(1)}} \cdots \partial {x}_{j_{\gamma_1}^{(1)}} } }_{\gamma_1 \text{times}} \cdots \underbrace{\frac{\partial^{\gamma_d}}{ \partial {x}_{j_1^{(d)}} \cdots \partial {x}_{j_{\gamma_d}^{(d)}} }}_{\gamma_d \text{times}} v^{}_{\eta} \Biggr| \\
&\hspace{-1.5cm}  \leq c | \det ({A}_T) | \| \widetilde{{A}}^{-1} \|_2 \alpha_k^{-1} \sum_{\eta=1}^d \sum_{|\varepsilon| = |\beta|+|\gamma|} \mathscr{H}^{\varepsilon} |\partial^{\varepsilon}_{x} v_{\eta}|.
\end{align*}
Because $1 \leq p \< \infty$, it holds that, for $1 \leq k \leq d$,
\begin{align*}
\displaystyle
\left \|  \partial^{\beta}_{\hat{x}} \partial^{\gamma}_{\hat{x}}  \hat{v}_k  \right \|_{L^{p}(\widehat{T})}^p 
&\leq c  | \det ({A}_T) |^{p-1} \| \widetilde{{A}}^{-1} \|_2^p \alpha_k^{-p} \sum_{|\varepsilon| = |\beta|+|\gamma|} \mathscr{H}^{\varepsilon p} \int_{T} |\partial^{\varepsilon}_{x} v |^p dx,
\end{align*}
which leads to \eqref{RT43} together with \eqref{jensen}, $T = \Phi_{T_0}^{-1}(T_0)$ and $v = \Psi_{T_0}^{-1} v_0$.
\qed	
\end{pf*}

\begin{Rem} \label{ex=01}
In inequality \eqref{RT43}, it is possible to obtain the estimates in $T_0$ by explicitly determining the matrix ${A}_{T_0}$.

Let $\hat{v} \in \mathcal{C}^{1}(\widehat{T})^d$ with ${v} = {\Psi}_T \hat{v}$ and  ${v}_0 = {\Psi}_{T_0} {v}$. Using \eqref{Anorm}, \eqref{CN331c}, \eqref{CN332} and the definition of Piola transformations, we have, for $1 \leq i,k \leq d$,
\begin{align*}
\displaystyle
\left| \frac{\partial \hat{v}_k}{\partial \hat{x}_{i}}\right|
 &\leq c  |\det ({A}_T) | \| \widetilde{{A}}^{-1} \|_2 \alpha_k^{-1} \sum_{\nu=1}^d \sum_{i_1^{(1)}, i_1^{(0,1)} = 1}^d \mathscr{H}_{i_1^{(1)}} | [{A}_{T_0}]_{{i_1^{(0,1)}} {i_1^{(1)}}} |  \left| \frac{\partial v_{0,\nu}}{\partial x^{(0)}_{{i_1^{(0,1)}}}} \right |.
 \end{align*}
Let $d=3$. We define the matrix ${A}_{T_0}$ as 
\begin{align*}
\displaystyle
{A}_{T_0} := 
\begin{pmatrix}
\cos \frac{\pi}{2}  & - \sin \frac{\pi}{2} & 0\\
 \sin \frac{\pi}{2} & \cos \frac{\pi}{2} & 0 \\
 0 & 0 & 1
\end{pmatrix}.
\end{align*}
We then have
\begin{align*}
\displaystyle
\left| \frac{\partial \hat{v}_k}{\partial \hat{x}_{i}}\right|
 &\leq c  |\det ({A}_T) | \| \widetilde{{A}}^{-1} \|_2 \alpha_k^{-1} \sum_{\nu=1}^3 \left( \mathscr{H}_{1}  \left| \frac{\partial v_{0,\nu}}{\partial x^{(0)}_{2}} \right | +  \mathscr{H}_{2}  \left| \frac{\partial v_{0,\nu}}{\partial x^{(0)}_{1}} \right | +  \mathscr{H}_{3}  \left| \frac{\partial v_{0,\nu}}{\partial x^{(0)}_{3}} \right | \right).
 \end{align*}
Because $1 \leq p \< \infty$, it holds that, for $1 \leq i, k \leq 3$,
\begin{align*}
\displaystyle
\left \| \frac{\partial \hat{v}_k}{\partial \hat{x}_{i}}\right \|^p_{L^p(\widehat{T})}
&\leq c  |\det ({A}_T) |^{p-1} \| \widetilde{{A}}^{-1} \|^p_2 \alpha_k^{-p} \\
&\hspace{-0.5cm} \times \left( \mathscr{H}_1^p \left \| \frac{\partial v_0}{\partial x^{(0)}_2} \right \|^p_{L^p(T)} + \mathscr{H}_2^p \left \| \frac{\partial v_0}{\partial x^{(0)}_1} \right \|^p_{L^p(T)} + \mathscr{H}_3^p \left \| \frac{\partial v_0}{\partial x^{(0)}_3} \right \|^p_{L^p(T)} \right).
 \end{align*}
\end{Rem}

The following two lemmata are divided into the element on $\mathfrak{T}^{(2)}$ or $\mathfrak{T}_1^{(3)}$ and the element on $\mathfrak{T}_2^{(3)}$. 

\begin{lem} \label{lem1142}
Let $T \in \mathfrak{T}^{(2)}$ or $T \in \mathfrak{T}_1^{(3)}$ satisfy Condition \ref{cond1} or Condition \ref{cond2}, respectively. Let $T_0 \subset \mathbb{R}^d$ be a simplex such that $T = \Phi_{T_0}^{-1}(T_0)$. Let $\beta := (\beta_1,\ldots,\beta_d) \in \mathbb{N}_0^d$ be a multi-index with $|\beta| = \ell$. Let $p \in [1,\infty)$. It then holds that, for any $\hat{v} = (\hat{v}_1,\ldots,\hat{v}_d)^T \in W^{\ell+1,p}(\widehat{T})^d$ with ${v} = ({v}_1,\ldots,{v}_d)^T := {\Psi}_{T} \hat{v}$ and ${v}_0 = (v_{0,1}, \ldots,v_{0,d})^T := \Psi_{T_0}{v}$,  
\begin{align}
\displaystyle
\left \| \partial^{\beta}_{\hat{x}} \nabla_{\hat{x}} \cdot \hat{v} \right \|_{L^p(\widehat{T})}
&\leq c |\det ({A}_{T})|^{\frac{p-1}{p}} \sum_{|\varepsilon| = \ell} \alpha^{\varepsilon} \left \|  \partial_{r^{(0)}}^{\varepsilon} \nabla_{x^{(0)}} \cdot v_0 \right \|_{L^p(T_0)}. \label{RT12}
\end{align}
If Assumption \ref{ass1} is imposed, it holds that
\begin{align}
\displaystyle
\left \| \partial^{\beta}_{\hat{x}}  \nabla_{\hat{x}} \cdot \hat{v}  \right \|_{L^p(\widehat{T})}
&\leq c |\det ({A}_T)|^{\frac{p-1}{p}} \sum_{|\varepsilon| = \ell} \mathscr{H}^{\varepsilon} \left \|  \partial^{\varepsilon}_{x} \nabla_{x} \cdot (\Psi_{T_0}^{-1}  v_0) \right \|_{L^p(\Phi_{T_0}^{-1}(T_0))}. \label{RT13}
\end{align}
\end{lem}

\begin{pf*}
Because the space $\mathcal{C}^{\ell+1}(\widehat{T})^d$ is dense
in the space $W^{\ell+1,p}(\widehat{T})^d$,  we show \eqref{RT12} and \eqref{RT13} for $\hat{v} \in \mathcal{C}^{\ell+1}(\widehat{T})^d$ with ${v} = {\Psi}_{T} \hat{v}$ and ${v}_0 = \Psi_{T_0}{v}$.

For a general derivative ${\partial}_{\hat{x}}^{\beta} \nabla_{\hat{x}} \cdot \hat{v}$ with order $|\beta| = \ell$, we obtaint, using \eqref{CN332},
 \begin{align*}
\displaystyle
| \partial^{\beta}_{\hat{x}} \nabla_{\hat{x}} \cdot \hat{v} | &=  \left| \frac{\partial^{|\beta|}}{\partial \hat{x}_1^{\beta_1} \cdots \partial \hat{x}_d^{\beta_d}} \nabla_{\hat{x}} \cdot \hat{v} \right| \notag\\
&\hspace{-1.5cm}  \leq | \det ({A}_T) | | \det ({A}_{T_0}) | \notag\\
&\hspace{-1.2cm} \Biggl | \underbrace{\sum_{i_1^{(1)},i_1^{(0,1)} = 1}^d \alpha_1 [{A}_{T_0}]_{i_1^{(0,1)} i_1^{(1)}} (r_1)_{i_1^{(1)}} \cdots   \sum_{i_{\beta_1}^{(1)},i_{\beta_1}^{(0,1)} = 1}^d \alpha_1 [{A}_{T_0}]_{i_{\beta_1}^{(0,1)} i_{\beta_1}^{(1)}} (r_1)_{i_{\beta_1}^{(1)}}  }_{\beta_1 \text{times}} \cdots \notag \\
&\hspace{-1.2cm} \underbrace{ \sum_{i_1^{(d)} , i_1^{(0,d)} = 1}^d \alpha_d [{A}_{T_0}]_{i_1^{(0,d)} i_1^{(d)}} (r_d)_{i_1^{(d)}}  \cdots \sum_{i_{\beta_d}^{(d)} , i_{\beta_d}^{(0,d)}= 1}^d \alpha_d [{A}_{T_0}]_{i_{\beta_d}^{(0,d)} i_{\beta_d}^{(d)}} (r_d)_{i_{\beta_d}^{(d)}} }_{\beta_d \text{times}}  \notag \\
&\hspace{-1.2cm}  \underbrace{\frac{\partial^{\beta_1}}{\partial {x}_{i_1^{(0,1)}}^{(0)} \cdots \partial {x}_{i_{\beta_1}^{(0,1)}}^{(0)}}}_{\beta_1 \text{times}} \cdots \underbrace{\frac{\partial^{\beta_d}}{\partial {x}_{i_1^{(0,d)}}^{(0)} \cdots \partial {x}_{i_{\beta_d}^{(0,d)}}^{(0)}}}_{\beta_d \text{times}} \nabla_{x^{(0)}} \cdot v_0 \Biggr | \\
&\hspace{-1.5cm} \leq |\det ({A}_T)| \sum_{|\varepsilon| = \ell} \alpha^{\varepsilon} |\partial_{r^{(0)}}^{\varepsilon} \nabla_{x^{(0)}} \cdot v_0 |
\end{align*}
It then holds that, using \eqref{jensen},
\begin{align*}
\displaystyle
\| \partial^{\beta}_{\hat{x}} \nabla_{\hat{x}} \cdot \hat{v} \|_{L^p(\widehat{T})}
&\leq c |\det ({A}_T)|^{\frac{p-1}{p}}\sum_{|\varepsilon| = \ell} \alpha^{\varepsilon} \| \partial_{r^{(0)}}^{\varepsilon} \nabla_{x^{(0)}} \cdot v_0 \|_{L^p(T_0)},
\end{align*}
which is the desired inequelity.

By an analogous argument, if Assumption \ref{ass1} is imposed, for a general derivative ${\partial}_{\hat{x}}^{\beta} \nabla_{\hat{x}} \cdot \hat{v}$ with order $|\beta| = \ell$, we have, using Note \ref{defi=theta}, 
\begin{align*}
\displaystyle
| \partial^{\beta}_{\hat{x}} \nabla_{\hat{x}} \cdot \hat{v} | &= \left| \frac{\partial^{|\beta|}}{\partial \hat{x}_1^{\beta_1} \cdots \partial \hat{x}_d^{\beta_d}} \nabla_{\hat{x}} \cdot \hat{v} \right| \notag\\
&\hspace{-1.5cm}  \leq c | \det ({A}_T) | \notag \\
&\hspace{-1.2cm}  \underbrace{\sum_{i_1^{(1)} = 1}^d \alpha_1 | [{\mathcal{A}}]_{i_1^{(1)} 1} | \cdots   \sum_{i_{\beta_1}^{(1)} = 1}^d \alpha_1 | [\widetilde{{A}}]_{i_{\beta_1}^{(1)} 1} |  }_{\beta_1 \text{times}} \cdots \underbrace{ \sum_{i_1^{(d)}  = 1}^d \alpha_d | [\widetilde{{A}}]_{i_1^{(d)} d} | \cdots \sum_{i_{\beta_d}^{(d)} = 1}^d \alpha_d | [\widetilde{{A}}]_{i_{\beta_d}^{(d)} d} | }_{\beta_d \text{times}}  \notag \\
&\hspace{-1.2cm} \Biggl | \underbrace{ \frac{\partial^{\beta_1}}{\partial {x}_{i_1^{(1)}} \cdots \partial {x}_{i_{\beta_1}^{(1)}}}}_{\beta_1 \text{times}} \cdots \underbrace{\frac{\partial^{\beta_d}}{\partial {x}_{i_1^{(d)}} \cdots \partial {x}_{i_{\beta_d}^{(d)}}}}_{\beta_d \text{times}} \nabla_{x} \cdot v \Biggr | \\
&\hspace{-1.5cm} \leq c |\det ({A}_T)| \sum_{|\varepsilon| = \ell} \mathscr{H}^{\varepsilon} | \partial_{x}^{\varepsilon}  \nabla_{x} \cdot v |.
\end{align*}
It then holds that, using \eqref{CN332} and \eqref{jensen},
\begin{align*}
\displaystyle
\| \partial^{\beta}_{\hat{x}} \nabla_{\hat{x}} \cdot \hat{v} \|_{L^p(\widehat{T})}
&\leq c |\det ({A}_T)|^{\frac{p-1}{p}}\sum_{|\varepsilon| = \ell} \mathscr{H}^{\varepsilon} \| \partial_{x}^{\varepsilon} \nabla_{x} \cdot v \|_{L^p(T)},
\end{align*}
which leads to \eqref{RT13} together with $T = \Phi_{T_0}^{-1}(T_0)$ and $v = \Psi_{T_0}^{-1} v_0$.
\qed	
\end{pf*}

\begin{lem} \label{lem1143}
Let $d=3$. Let $T \in \mathfrak{T}_2^{(3)}$ satisfy  Condition \ref{cond2}. Let $T_0 \subset \mathbb{R}^3$ be a simplex such that $T = \Phi_{T_0}^{-1}(T_0)$. Let $\ell \in \mathbb{N}_0$ and  $k \in \mathbb{N}$ with $1 \leq k \leq 3$. Let $\beta := (\beta_1,\beta_2,\beta_3) \in \mathbb{N}_0^3$ be a multi-index with $|\beta| = \ell$. Let $p \in [1,\infty)$. It holds that, for any $\hat{v} = (\hat{v}_1,\hat{v}_2,\hat{v}_3)^T \in W^{\ell+1,p}(\widehat{T})^3$ with ${v} = ({v}_1,v_2,{v}_3)^T := {\Psi}_{T} \hat{v}$ and $v_0 = (v_{0,1},v_{0,2},v_{0,3})^T := \Psi_{T_0} v$, 
\begin{align}
\displaystyle
\left \| \partial_{\hat{x}}^{\beta} \frac{\partial \hat{v}_k}{\partial \hat{x}_k} \right\|_{L^p(\widehat{T})} 
&\leq c  |\det ({A}_T)|^{\frac{p-1}{p}} \| \widetilde{{A}}^{-1} \|_2 \sum_{|\varepsilon| = \ell} \alpha^{\varepsilon} \left \| \partial^{\varepsilon}_{r^{(0)}} \frac{\partial v_0}{\partial r_k^{(0)}} \right \|_{L^p(T_0)^3}.  \label{RT14}
\end{align}
If Assumption \ref{ass1} is imposed, it holds that
\begin{align}
\displaystyle
\left \| \partial^{\beta}_{\hat{x}}  \frac{\partial \hat{v}_k}{\partial \hat{x}_{k}} \right \|_{L^p(\widehat{T})}
&\leq c  |\det ({A}_T)|^{\frac{p-1}{p}} \| \widetilde{{A}}^{-1} \|_2 \sum_{|\varepsilon| = \ell} \mathscr{H}^{\varepsilon} \left \| \partial^{\varepsilon}_{x} \frac{\partial (\Psi_{T_0}^{-1} v)}{\partial r_k} \right \|_{L^p(\Phi_{T_0}^{-1}(T_0))^3}. \label{RT14b}
\end{align}
\end{lem}

\begin{pf*}
Because the space $\mathcal{C}^{\ell+1}(\widehat{T})^3$ is dense
in the space $W^{\ell+1,p}(\widehat{T})^3$,  we show \eqref{RT12} and \eqref{RT13} for $\hat{v} \in \mathcal{C}^{\ell+1}(\widehat{T})^3$ with ${v} = {\Psi}_{T} \hat{v}$ and ${v}_0 = \Psi_{T_0}{v}$.

For a general derivative ${\partial}_{\hat{x}}^{\beta} \frac{\partial \hat{v}_k}{\partial \hat{x}_{k}}$ ($1 \leq k \leq 3$) with order $|\beta| = \ell$, we obtain, using \eqref{Anorm}, \eqref{CN331c} and \eqref{CN332},
\begin{align*}
\displaystyle
\left| \partial^{\beta}_{\hat{x}}  \frac{\partial \hat{v}_k}{\partial \hat{x}_{k}} \right| &= \left|  \frac{\partial^{|\beta| }}{\partial \hat{x}_1^{\beta_1} \partial \hat{x}_2^{\beta_2} \partial \hat{x}_3^{\beta_3} }   \frac{\partial \hat{v}_k}{\partial \hat{x}_{k}} \right| \notag\\
&\hspace{-1.5cm}  \leq | \det ({A}_T) |  | \det ({A}_{T_0}) |  \sum_{\eta,\nu=1}^3 | [ \widetilde{{A}}^{-1}]_{k \eta} |  | [ {A}_{T_0}^{-1}]_{\eta \nu} | \notag\\
&\hspace{-1.2cm} \Biggl | \underbrace{\sum_{i_1^{(1)},i_1^{(0,1)} = 1}^3 \alpha_1  [{A}_{T_0}]_{i_1^{(0,1)} i_1^{(1)}} (r_1)_{i_1^{(1)}} \cdots   \sum_{i_{\beta_1}^{(1)},i_{\beta_1}^{(0,1)} = 1}^3 \alpha_1 [{A}_{T_0}]_{i_{\beta_1}^{(0,1)} i_{\beta_1}^{(1)}} (r_1)_{i_{\beta_1}^{(1)}}   }_{\beta_1 \text{times}} \notag \\
&\hspace{-1.2cm}  \underbrace{\sum_{i_1^{(2)},i_1^{(0,2)} = 1}^3 \alpha_2 [{A}_{T_0}]_{i_1^{(0,2)} i_1^{(2)}} (r_2)_{i_1^{(2)}} \cdots   \sum_{i_{\beta_2}^{(2)},i_{\beta_2}^{(0,2)} = 1}^3 \alpha_2 [{A}_{T_0}]_{i_{\beta_2}^{(0,2)} i_{\beta_2}^{(2)}} (r_2)_{i_{\beta_2}^{(2)}}   }_{\beta_2 \text{times}}  \notag \\
&\hspace{-1.2cm} \underbrace{ \sum_{i_1^{(3)} , i_1^{(0,3)} = 1}^3 \alpha_3  [{A}_{T_0}]_{i_1^{(0,3)} i_1^{(3)}} (r_3)_{i_1^{(3)}}  \cdots \sum_{i_{\beta_3}^{(3)} , i_{\beta_3}^{(0,3)}= 1}^3 \alpha_3 [{A}_{T_0}]_{i_{\beta_3}^{(0,3)} i_{\beta_3}^{(3)}} (r_3)_{i_{\beta_3}^{(3)}} }_{\beta_3 \text{times}}  \notag \\
&\hspace{-1.2cm}  \underbrace{\frac{\partial^{\beta_1}}{\partial {x}_{i_1^{(0,1)}}^{(0)} \cdots \partial {x}_{i_{\beta_1}^{(0,1)}}^{(0)}}}_{\beta_1 \text{times}} \underbrace{\frac{\partial^{\beta_2}}{\partial {x}_{i_1^{(0,2)}}^{(0)} \cdots \partial {x}_{i_{\beta_2}^{(0,2)}}^{(0)}}}_{\beta_2 \text{times}} \underbrace{\frac{\partial^{\beta_3}}{\partial {x}_{i_1^{(0,3)}}^{(0)} \cdots \partial {x}_{i_{\beta_3}^{(0,3)}}^{(0)}}}_{\beta_3 \text{times}}   \frac{\partial v_{0,\nu}}{\partial r_k^{(0)}}  \Biggr | \\
&\hspace{-1.5cm}  \leq c | \det ({A}_T) | \| \widetilde{\mathcal{A}}^{-1} \|_2  \sum_{\nu=1}^3 \sum_{|\varepsilon| = |\beta|} \alpha^{\varepsilon} \left| \partial^{\varepsilon}_{r^{(0)}} \frac{\partial v_{0,\nu}}{\partial r_k^{r^{(0)}}} \right|.
\end{align*}
It then holds that, using \eqref{jensen},
\begin{align*}
\displaystyle
\left \| \partial^{\beta}_{\hat{x}}  \frac{\partial \hat{v}_k}{\partial \hat{x}_{k}} \right \|_{L^p(\widehat{T})}
&\leq c  |\det ({A}_T)|^{\frac{p-1}{p}} \| \widetilde{{A}}^{-1} \|_2 \sum_{|\varepsilon| = \ell} \alpha^{\varepsilon} \left \| \partial^{\varepsilon}_{r^{(0)}} \frac{\partial v_0}{\partial r_k^{r^{(0)}}} \right \|_{L^p(T_0)^3},
\end{align*}
which is the desired inequelity.

By an analogous argument, if Assumption \ref{ass1} is imposed, for a general derivative ${\partial}_{\hat{x}}^{\beta} \frac{\partial \hat{v}_k}{\partial \hat{x}_{k}}$ ($1 \leq k \leq 3$) with order $|\beta| = \ell$, we have, using Note \ref{defi=theta} and \eqref{Anorm}, 
\begin{align*}
\displaystyle
\left| \partial^{\beta}_{\hat{x}}  \frac{\partial \hat{v}_k}{\partial \hat{x}_{k}} \right| &= \left | \frac{\partial^{|\beta| }}{\partial \hat{x}_1^{\beta_1} \partial \hat{x}_2^{\beta_2} \partial \hat{x}_3^{\beta_3} }   \frac{\partial \hat{v}_k}{\partial \hat{x}_{k}} \right| \notag\\
&\hspace{-1.5cm} \leq | \det ({A}_T) | \sum_{\eta=1}^3 | [ \widetilde{{A}}^{-1}]_{k \eta} |   \notag \\
&\hspace{-1.2cm}  \underbrace{\sum_{i_1^{(1)} = 1}^3 \alpha_1 | [\widetilde{{A}}]_{i_1^{(1)} 1}| \cdots   \sum_{i_{\beta_1}^{(1)} = 1}^3 \alpha_1 | [\widetilde{{A}}]_{i_{\beta_1}^{(1)} 1} |  }_{\beta_1 \text{times}}\underbrace{ \sum_{i_1^{(2)}  = 1}^3 \alpha_2 | [\widetilde{{A}}]_{i_1^{(2)} 2} | \cdots \sum_{i_{\beta_2}^{(2)} = 1}^3 \alpha_2 | [\widetilde{{A}}]_{i_{\beta_2}^{(2)} 2} | }_{\beta_2 \text{times}}  \notag \\
&\hspace{-1.2cm} \underbrace{ \sum_{i_1^{(3)}  = 1}^3 \alpha_3 | [\widetilde{{A}}]_{i_1^{(3)} 3} | \cdots \sum_{i_{\beta_3}^{(3)} = 1}^3 \alpha_3 | [\widetilde{{A}}]_{i_{\beta_3}^{(3)} 3} | }_{\beta_3 \text{times}}  \notag \\
&\hspace{-1.2cm} \Biggl | \underbrace{\frac{\partial^{\beta_1}}{\partial {x}_{i_1^{(1)}} \cdots \partial {x}_{i_{\beta_1}^{(1)}}}}_{\beta_1 \text{times}}  \underbrace{\frac{\partial^{\beta_2}}{\partial {x}_{i_1^{(2)}} \cdots \partial {x}_{i_{\beta_2}^{(2)}}}}_{\beta_2 \text{times}} \underbrace{\frac{\partial^{\beta_3}}{\partial {x}_{i_1^{(3)}} \cdots \partial {x}_{i_{\beta_3}^{(3)}}}}_{\beta_3 \text{times}} \frac{\partial v_{\eta}}{ \partial r_{k}} \Biggr| \\
&\hspace{-1.5cm} \leq c | \det ({A}_T) | \| \widetilde{{A}}^{-1} \|_2 \sum_{\eta=1}^3 \sum_{|\varepsilon| = |\beta|} \mathscr{H}^{\varepsilon} \left | \partial_{x}^{\varepsilon} \frac{\partial v_{\eta}}{ \partial r_{k}} \right|.
\end{align*}
It then holds that, using \eqref{jensen},
\begin{align*}
\displaystyle
\left \| \partial^{\beta}_{\hat{x}}  \frac{\partial \hat{v}_k}{\partial \hat{x}_{k}} \right \|_{L^p(\widehat{T})}
&\leq c  |\det ({A}_T)|^{\frac{p-1}{p}} \| \widetilde{{A}}^{-1} \|_2 \sum_{|\varepsilon| = \ell} \mathscr{H}^{\varepsilon} \left \| \partial^{\varepsilon}_{x} \frac{\partial v}{\partial r_k} \right \|_{L^p(T)^3},
\end{align*}
which leads to \eqref{RT14b} together with $T = \Phi_{T_0}^{-1}(T_0)$ and $v = \Psi_{T_0}^{-1} v_0$.
\qed	
\end{pf*}

\section{Local Raviart--Thomas interpolation error estimates}
We introduce component-wise stability for the Raviart--Thomas interpolation proposed in \cite{AcoApe10}; see also \cite{AcoDur99}. 

We first introduce component-wise stability estimates in the reference element $\widehat{T}_1 = \conv \{ 0,e_1, \ldots,e_d \}$, where $e_1, \ldots, e_d \in \mathbb{R}^d$ are the canonical basis.

\begin{lem} \label{rt=lem7}
For $k \in \mathbb{N}_0$, there exists a constant $C_1^{(i)}(k)$, $i=1,\dots,d$ such that for all $\hat{u} = (\hat{u}_1, \ldots, \hat{u}_d)^T \in W^{1,p}(\widehat{T}_1)^d$, 
\begin{align}
\displaystyle
\| (I_{\widehat{T}_1}^{RT^k} \hat{u})_i \|_{L^p(\widehat{T}_1)} \leq C_1^{(i)}(k) \left( \| \hat{u}_i \|_{W^{1,p}(\widehat{T}_1)} + \| \nabla_{\hat{x}} \cdot \hat{u} \|_{L^p(\widehat{T}_1)} \right), \quad i=1, \ldots, d.\label{rt4}
\end{align}
\end{lem}

\begin{pf*}
The proof is provided in \cite[Lemma 3.3]{AcoApe10} for the case $d=3$. The estimate in the case $d=2$ can be proved analogously.
\qed
\end{pf*}

We next provide component-wise stability estimates in the reference element $\widehat{T}_2 = \conv \{ 0,e_1, e_1 + e_2 , e_3 \}$. 

\begin{lem} \label{rt=lem9}
For $k \in \mathbb{N}_0$, there exists a constant $C_2^{(i)}(k)$, $i=1,2,3$ such that, for all $\hat{u} = (\hat{u}_1,\hat{u}_2, \hat{u}_3)^T \in W^{1,p}(\widehat{T}_2)^3$, 
\begin{align}
\displaystyle
\| (I_{\widehat{T}_2}^{RT^k} \hat{u})_i \|_{L^p(\widehat{T}_2)} 
&\leq C_2^{(i)}(k) \left( \| \hat{u}_i \|_{W^{1,p}(\widehat{T}_2)} + \sum_{j=1, j\neq i}^3 \left\| \frac{\partial \hat{u}_j}{\partial \hat{x}_j} \right\|_{L^p(\widehat{T}_2)} \right), \quad i=1,2,3. \label{rt5}
\end{align}
\end{lem}

\begin{pf*}
The proof is provided in \cite[Lemma 4.3]{AcoApe10}. We remark that our reference element, in this case, is different from that in  \cite[Lemma 4.3]{AcoApe10}. However, the lemma can be proved using an analogous argument.
\qed
\end{pf*}

\subsection{Stability of the local Raviart--Thomas interpolation}

%%三角不等式を用いて誤差評価から導ける

\begin{lem} \label{lem1151}
Let $p \in [1,\infty)$. Let $T \in \mathfrak{T}^{(2)}$ or $T \in \mathfrak{T}_1^{(3)}$ satisfy Condition \ref{cond1} or Condition \ref{cond2}, respectively. Let $T_0 \subset \mathbb{R}^d$ be a simplex such that $T = \Phi_{T_0}^{-1}(T_0)$. It holds that, for any $\hat{v} = (\hat{v}_1,\ldots,\hat{v}_d)^T \in W^{1,p}(\widehat{T})^d$ with ${v} = ({v}_1,\ldots,{v}_d)^T := {\Psi}_{T} \hat{v}$ and ${v}_0 = (v_{0,1}, \ldots,v_{0,d})^T := \Psi_{T_0}{v}$, 
\begin{align}
\displaystyle
&\| {I_{T_0}^{RT^k} v_0} \|_{L^p({T_0})^d} \notag \\
&\quad \leq c \left[ \frac{H_{T_0}}{h_{T_0}} \left( \| v_0 \|_{L^p(T_0)^d} + \sum_{|\varepsilon|=1} \alpha^{\varepsilon} \left \| \partial_{r^{(0)}}^{\varepsilon} v_0 \right \|_{L^p(T_0)^d} \right ) + h_{T_0} \| \nabla_{x^{(0)}} \cdot {v}_0 \|_{L^p({T}_0)} \right]. \label{RT51}
\end{align}

Let $d=3$. Let $T \in \mathfrak{T}_2^{(3)}$ satisfy Condition \ref{cond2}. Let $T_0 \subset \mathbb{R}^3$ be a simplex such that $T = \Phi_{T_0}^{-1}(T_0)$.  It holds that, for any $\hat{v} = (\hat{v}_1,\hat{v}_2,\hat{v}_3)^T \in W^{1,p}(\widehat{T})^3$ with ${v} = ({v}_1,v_2,{v}_3)^T := {\Psi}_{T} \hat{v}$ and ${v}_0 = (v_{0,1}, v_{0,2},v_{0,3})^T := \Psi_{T_0}{v}$, 
\begin{align}
\displaystyle
\| {I_{T_0}^{RT^k} v_0} \|_{L^p({T_0})^3} 
&\leq c \frac{H_{T_0}}{h_{T_0}} \Biggl[  \| v_0 \|_{L^p(T_0)^3} + h_{T_0} \sum_{k=1}^3 \left \| \frac{\partial v_0}{\partial r_k^{(0)}} \right \|_{L^p(T_0)^3} \Biggr]. \label{RT58}
\end{align}
\end{lem}

\begin{pf*}
The inequality \eqref{RT51} follows from \eqref{RT41}, the component-wise stability \eqref{rt4}, \eqref{RT42} with $\ell=0$ and $m \in \{0,1 \}$, and \eqref{RT12} with $\ell=0$. The inequality \eqref{RT58} follows from \eqref{RT41}, the component-wise stability \eqref{rt5}, \eqref{RT42} with $\ell=0$ and $m \in \{0,1 \}$, and \eqref{RT14} with $\ell=0$.
\qed	
\end{pf*}

\subsection{Remarks on anisotropic interpolation error analysis} \label{Sec62}
In the proof of Theorem 3 in \cite{IshKobTsu21a}, we used  the following estimate in \cite[Lemmas 6 and 7]{IshKobTsu21a}: For any $v \in H^1(T)^d$,
\begin{align*}
\displaystyle
\frac{\| I_T^{RT} v - v \|_{L^2({T})^d}}{| {v} |_{H^{1}({T})^d}} 
&\leq C^{P,d} \frac{H_T}{h_T} h_T \frac{\left( \sum_{i=1}^d \alpha_i^2 \|  ( I_{\widehat{T}}^{RT}  \hat{v} )_i - \hat{v}_i  \|_{L^2(\widehat{T})}^2 \right)^{1/2}}{\left( \sum_{i=1}^d \alpha_i^2 |\hat{v}_i |_{H^{1}(\widehat{T})}^2 \right)^{1/2}}.
\end{align*}
If the component-wise stability of the Raviart--Thomas interpolation on the reference element $\widehat{T}$
\begin{align}
\displaystyle
 \| ( I_{\widehat{T}}^{RT}  \hat{v} )_i - \hat{v}_i  \|_{L^2(\widehat{T})} \leq c | \hat{v}_i|_{H^{1}(\widehat{T})}, \quad i=1,\ldots,d \label{intro1}
\end{align}
holds, then the target estimate
\begin{align*}
\displaystyle
\| I_T^{RT} v - v \|_{L^2({T})^d}
\leq c \frac{H_T}{h_T} h_T | {v} |_{H^{1}({T})^d}
\end{align*}
holds. However, the estimate \eqref{intro1} generally does not hold (see \cite[Introduction]{AcoApe10}): that is, we cannot apply the Babu\v{s}ka and Aziz technique \cite{BabAzi76}. We provide a counterexample of  \cite[Introduction]{AcoApe10}. 

We consider the simplex $\widehat{T} \subset \mathbb{R}^2$ with vertices $\widehat{P}_1 := (0,0)^T$, $\widehat{P}_2 := (1,0)^T$, and $\widehat{P}_3 := (0,1)^T$. For $1 \leq i \leq 3$, let $\widehat{F}_i$ be the face of $\widehat{T}$ opposite to $\widehat{P}_i$. The Raviart--Thomas interpolation of $\hat{v}$ is defined as
\begin{align*}
\displaystyle
 I_{\widehat{T}}^{RT} \hat{v}
  = \sum_{i=1}^3 \left( \int_{\widehat{F}_i} \hat{v} \cdot \hat{n}_i d \hat{s} \right) \hat{\theta}_i \in RT^0,
% = \sum_{i=1}^3 \frac{1}{|\widehat{F}_i|} \left( \int_{\widehat{F}_i} \hat{v} \cdot \hat{n}_i d \hat{s} \right) \hat{\theta}_i \in RT^0,
\end{align*}
where 
\begin{align*}
\displaystyle
\hat{\theta}_i := \frac{1}{2 |\widehat{T}|} (\hat{x} - \widehat{P}_i), \quad \hat{x} = (\hat{x}_1,\hat{x}_2)^T.
\end{align*}
Setting $\hat{v} := (0, \hat{x}_2^2)^T$  yields $ I_{\widehat{T}}^{RT} \hat{v} = \frac{1}{3} (\hat{x}_1,\hat{x}_2)^T$.
%\begin{align*}
%\displaystyle
 %I_{\widehat{T}}^{RT} \hat{v}
 %&= \frac{1}{\sqrt{2}} \left( \int_{\widehat{F}_1} \hat{x}_1^2 d \hat{s} \right) (\hat{x}_1,\hat{x}_2)^T - \left( \int_{\widehat{F}_3} \hat{x}_1^2 d \hat{s} \right)  (\hat{x}_1,\hat{x}_2 - 1)^T \\
% &= \frac{1}{3} (\hat{x}_1,\hat{x}_2)^T.
%\end{align*}
%\begin{align*}
%\displaystyle
% I_{\widehat{T}}^{RT} ( \hat{v} - \hat{q} )
 %&= \frac{\sqrt{2}}{12} (\hat{x}_1 , \hat{x}_2)^T - \frac{1}{3} (\hat{x}_1 , \hat{x}_2 - 1)^T.
%\end{align*}
This implies that $( I_{\widehat{T}}^{RT} \hat{v} )_1 - \hat{v}_1 \not\equiv 0$ for any $\hat{x} \in \mathbb{R}^2$.

%As described in Section 1, to overcome this difficulty, we use the component-wise stability of the  Raviart--Thomas interpolation on reference elements in \cite{AcoDur99,AcoApe10}.

\subsection{Main theorems}
The Bramble--Hilbert--type lemma (e.g., see \cite{DupSco80,BreSco08}) plays a major role in interpolation error analysis. We use the following estimates on anisotropic meshes proposed in \cite[Lemma 2.1]{Ape99}.

\begin{lem} \label{lem9}
%Let $D := \cup_{j=1}^J D_j \subset \mathbb{R}^d$ with $d \in \{ 2,3\}$, be a connected open set that is the union of a finite collection of domains $D_j \subset \mathbb{R}^d$ that are star-shaped with respect to balls $B_j$. 
Let $D \subset \mathbb{R}^d$ with $d \in \{ 2,3\}$, be a connected open set that is star-shaped with respect to balls $B$. Let $\gamma$ be a multi-index with $m := |\gamma|$ and $\varphi \in L^1(D)$ be a function with $\partial^{\gamma} \varphi \in W^{\ell -m,p}(D)$, where $\ell \in \mathbb{N}$, $m \in \mathbb{N}_0$, $0 \leq m \leq \ell$, $p \in [1,\infty]$. It then holds that
\begin{align}
\displaystyle
\| \partial^{\gamma} (\varphi - Q^{(\ell)} \varphi) \|_{W^{\ell -m,p}(D)} \leq C^{BH} |\partial^{\gamma} \varphi|_{W^{\ell-m,p}(D)},  \label{lem9=1}
\end{align}
where $C^{BH}$ depends only on $d$, $\ell$, $\diam D$, and $\diam B$, and $ Q^{(\ell)} \varphi$ is defined as
\begin{align}
\displaystyle
(Q^{(\ell)} \varphi)(x) := \sum_{|\delta| \leq \ell -1} \int_B \eta(y) (\partial^{\delta}\varphi)(y) \frac{(x-y)^{\delta}}{\delta !} dy \in \mathcal{P}^{\ell -1},  \label{lem9=2}
\end{align}
where $\eta \in \mathcal{C}_0^{\infty}(B)$ is a given function with $\int_B \eta dx = 1$.
\end{lem}

The following two theorems are divided into the element on $\mathfrak{T}^{(2)}$ or $\mathfrak{T}_1^{(3)}$ and the element on $\mathfrak{T}_2^{(3)}$.

\begin{thr} \label{thr1161}
Let $p \in [1,\infty)$. Let $T \in \mathfrak{T}^{(2)}$ or $T \in \mathfrak{T}_1^{(3)}$ satisfy Condition \ref{cond1} or Condition \ref{cond2}, respectively.  Let $T_0 \subset \mathbb{R}^d$ be a simplex such that $T = \Phi_{T_0}^{-1}(T_0)$. For $k \in \mathbb{N}_0$, let $\{ {T}_0 , RT^k({T}_0) , {\Sigma}_0 \}$ be the Raviart--Thomas finite element and $I_{T_0}^{RT^k}$ the local interpolation operator defined in \eqref{RTlocal}. Let $\ell$ be such that $0 \leq \ell \leq k$. For any $\hat{v} \in W^{\ell+1,p}(\widehat{T})^d$ with ${v} = ({v}_1,\ldots,{v}_d)^T := {\Psi}_{T} \hat{v}$ and ${v}_0 = (v_{0,1}, \ldots,v_{0,d})^T := \Psi_{T_0}{v}$, it holds that
\begin{align}
\displaystyle
&\| I_{T_0}^{RT^k} v_0 - v_0 \|_{L^p(T_0)^d} \notag \\
&\leq  c \left( \frac{H_{T_0}}{h_{T_0}} \sum_{|\varepsilon| = \ell+ 1} \alpha^{\varepsilon} \left \| \partial_{r^{(0)}}^{\varepsilon} v_0 \right \|_{L^p(T_0)^d} +  h_{T_0} \sum_{|\beta| = \ell} \alpha^{\beta} \| \partial_{r^{(0)}}^{\beta} \nabla_{x^{(0)}} \cdot {v}_0 \|_{L^{p}({T}_0)} \right). \label{RT61}
\end{align}
If Assumption \ref{ass1} is imposed, it holds that
\begin{align}
\displaystyle
&\| I_{T_0}^{RT^k} v_0 - v_0 \|_{L^p(T_0)^d} \notag \\
&\leq c  \Biggl( \frac{H_{T_0}}{h_{T_0}}  \sum_{|\varepsilon| = \ell + 1} \mathscr{H}^{\varepsilon} \| \partial^{\varepsilon}_{x} (\Psi_{T_0}^{-1} v_0) \|_{L^p(\Phi_{T_1}^{-1}(T_0))^d} \notag \\
&\hspace{1cm} +  h_{T_0} \sum_{|\beta| = \ell} \mathscr{H}^{\beta} \| \partial^{\beta}_{x} \nabla_{x} \cdot (\Psi_{T_0}^{-1} v_0) \|_{L^{p}(\Phi_{T_0}^{-1}(T_0))} \Biggr). \label{RT62}
\end{align}
\end{thr}

\begin{pf*}
Let $\hat{v} \in W^{\ell+1,p}(\widehat{T})^d$. Let ${I}_{\widehat{T}}^{RT^k}$ be the local interpolation operators on $\widehat{T}_1$ defined by \eqref{RTinter1} and \eqref{RTinter2}. If ${q} \in \mathcal{P}^{\ell} ({T}_0)^d \subset RT^k({T}_0)$, then $I_{{T}_0}^{RT} {q} = {q}$. 

We set ${\mathfrak{Q}}^{(\ell+1)} {v_0} := ({Q}^{(\ell+1)} {v}_{0,1}, \ldots , {Q}^{(\ell+1)} {v}_{0,d})^T \in \mathcal{P}^{\ell}({T}_0)^d$, where ${Q}^{(\ell+1)} {v}_{0,j}$ is defined by \eqref{lem9=2} for any $j$. We then obtain
\begin{align}
\displaystyle
\| I_{T_0}^{RT^k} v_0 - v_0 \|_{L^p(T_0)^d} &\leq \| I_{T_0}^{RT^k} (v_0 - \mathfrak{Q}^{(\ell+1)} {v}_0) \|_{L^p(T_0)^d} + \| \mathfrak{Q}^{(\ell+1)} {v_0} - v_0 \|_{L^p(T_0)^d}. \label{RT63}
\end{align}
The inequality \eqref{RT41} for the first term on the right-hand side of \eqref{RT63} yield
\begin{align}
\displaystyle
&\| I_{T_0}^{RT^k} (v_0 - \mathfrak{Q}^{(\ell+1)} {v_0}) \|_{L^p(T_0)^d} \notag \\
&\quad \leq c  |\det({A}_T)|^{\frac{1-p }{p}} \| \widetilde{{A}} \|_{2} \left(  \sum_{j=1}^d \alpha_j^p \|  \{  I_{\widehat{T}}^{RT^k} (\hat{v} -  \widehat{\mathfrak{Q}}^{(\ell+1)} \hat{v}) \}_j \|_{L^p(\widehat{T})}^p \right)^{1/p}.\label{RT64}
\end{align}
The component-wise stability \eqref{rt4} and \eqref{jensen} yields
\begin{align}
\displaystyle
& \left(  \sum_{j=1}^d \alpha_j^p \|  \{  I_{\widehat{T}}^{RT^k} (\hat{v} -  \widehat{\mathfrak{Q}}^{(\ell+1)} \hat{v}) \}_j \|_{L^p(\widehat{T})}^p \right)^{1/p} \notag \\
&\quad \leq  \sum_{j=1}^d \alpha_j \|  \{  I_{\widehat{T}}^{RT^k} (\hat{v} -  \widehat{\mathfrak{Q}}^{(\ell+1)} \hat{v}) \}_j \|_{L^p(\widehat{T})} \notag \\
&\quad \leq c \sum_{j=1}^d \alpha_j \left( \| \hat{v}_j - \widehat{Q}^{(\ell+1)} \hat{v}_j \|_{W^{1,p}(\widehat{T})} + \| \nabla_{\hat{x}} \cdot ( \hat{v} - \widehat{\mathfrak{Q}}^{(\ell+1)} \hat{v}) \|_{L^p(\widehat{T})} \right). \label{RT65}
\end{align}
The inequalities \eqref{RT41} for the second term on the right-hand side of \eqref{RT63} and \eqref{jensen} yield
\begin{align}
\displaystyle
&\| \mathfrak{Q}^{(\ell+1)} {v}_0 - v_0 \|_{L^p(T_0)^d} \notag \\
&\quad \leq c  |\det({A}_T)|^{\frac{1-p }{p}} \| \widetilde{{A}} \|_{2} \left(  \sum_{j=1}^d \alpha_j^p \|  \widehat{Q}^{(\ell+1)} \hat{v}_j - \hat{v}_j \|_{L^p(\widehat{T})}^p \right)^{1/p} \notag \\
&\quad \leq c  |\det({A}_T)|^{\frac{1-p }{p}} \| \widetilde{{A}} \|_{2} \sum_{j=1}^d \alpha_j \|  \widehat{Q}^{(\ell+1)} \hat{v}_j - \hat{v}_j \|_{L^p(\widehat{T})}. \label{RT66}
\end{align}
The inequality \eqref{RT63} together with \eqref{RT64}, \eqref{RT65} and \eqref{RT66} leads to
\begin{align}
\displaystyle
&\| I_{T_0}^{RT^k} v_0 - v_0 \|_{L^p(T_0)^d} \notag \\
&\quad \leq |\det({A}_T)|^{\frac{1-p }{p}} \| \widetilde{{A}} \|_{2} \notag \\
&\quad \quad \sum_{j=1}^d \alpha_j \left( \| \hat{v}_j - \widehat{Q}^{(\ell+1)} \hat{v}_j \|_{W^{1,p}(\widehat{T})} + \| \nabla_{\hat{x}} \cdot ( \hat{v} - \widehat{\mathfrak{Q}}^{(\ell+1)} \hat{v}) \|_{L^p(\widehat{T})} \right). \label{RT63b}
\end{align}
The Bramble--Hilbert-type lemma (Lemma \ref{lem9}) and \eqref{RT42},
\begin{align}
\displaystyle
&\| \hat{v}_j - \widehat{Q}^{(\ell+1)} \hat{v}_j \|_{W^{1,p}(\widehat{T})}^p \notag\\
&\quad = \| \hat{v}_j - \widehat{Q}^{(\ell+1)} \hat{v}_j \|_{L^{p}(\widehat{T})}^p + \sum_{k=1}^d \left \| \frac{\partial}{\partial \hat{x}_k} ( \hat{v}_j - \widehat{Q}^{(\ell+1)} \hat{v}_j ) \right \|_{L^{p}(\widehat{T})}^p \notag\\
&\quad \leq c \left( \sum_{|\gamma| = \ell+ 1} \left \|  \partial_{\hat{x}}^{\gamma} \hat{v}_j \right \|_{L^{p}(\widehat{T})}^p + \sum_{k=1}^d \sum_{|\beta| = \ell} \left\|  \partial_{\hat{x}}^{\beta} \frac{\partial \hat{v}_j}{\partial \hat{x}_k} \right\|_{L^{p}(\widehat{T})}^p \right) \notag\\
&\quad \leq c  |\det ({A}_T)|^{p-1} \alpha_j^{-p} \| \widetilde{{A}}^{-1} \|_{2}^p \left( \sum_{|\varepsilon| = \ell+ 1} \alpha^{\varepsilon} \left \| \partial_{r^{(0)}}^{\varepsilon} v_0 \right \|_{L^p(T_0)^d} \right)^p. \label{RT67}
\end{align}
Because from \cite[Proposition 4.1.17]{BreSco08} it holds that
\begin{align}
\displaystyle
\nabla_{\hat{x}} \cdot ( \widehat{\mathfrak{Q}}^{(\ell+1)} \hat{v}) = \widehat{Q}^{\ell}(\nabla_{\hat{x}} \cdot \hat{v}), \label{RT68}
\end{align}
from Lemma \ref{lem9} and \eqref{RT12},
\begin{align}
\displaystyle
 \|  \nabla_{\hat{x}} \cdot ( \hat{v} - \widehat{\mathfrak{Q}}^{(\ell+1)} \hat{v})) \|_{L^p(\widehat{T})}^p
  &=  \| \nabla_{\hat{x}} \cdot \hat{v} -\widehat{Q}^{\ell}(\nabla_{\hat{x}} \cdot \hat{v}) \|_{L^p(\widehat{T})}^p \notag\\
 &\leq  \| \nabla_{\hat{x}} \cdot \hat{v} -\widehat{Q}^{\ell}( \nabla_{\hat{x}} \cdot \hat{v}) \|_{W^{\ell,p}(\widehat{T})}^p \notag\\
 &\leq c | \nabla_{\hat{x}} \cdot \hat{v} |_{W^{\ell,p}(\widehat{T})}^p 
 = c \sum_{|\beta| = \ell} \| \partial^{\beta} \nabla_{\hat{x}} \cdot \hat{v} \|_{L^{p}(\widehat{T})}^p \notag\\
& \leq  c  |\det ({A}_T)|^{{p-1}} \left( \sum_{|\varepsilon| = \ell} \alpha^{\varepsilon} \left \|  \partial_{r^{(0)}}^{\varepsilon} \nabla_{x^{(0)}} \cdot v_0 \right \|_{L^p(T_0)} \right)^p. \label{RT69}
\end{align}
Combining \eqref{RT63b}, \eqref{RT67}, and \eqref{RT69} with \eqref{CN331b} yields the target result \eqref{RT61}.

If Assumption \ref{ass1} is imposed, we use instead \eqref{RT43} of \eqref{RT42} and \eqref{RT13} instead of \eqref{RT12}. The inequality \eqref{RT62} then follows from \eqref{RT63b}, \eqref{RT43} and \eqref{RT13} with \eqref{CN331b}.
\qed	
\end{pf*}

\begin{thr} \label{thr1161b}
Let $d=3$ and $p \in [1,\infty)$. Let $T \in \mathfrak{T}_2^{(3)}$ satisfy Condition \ref{cond2}. Let $T_0 \subset \mathbb{R}^3$ be a simplex such that $T = \Phi_{T_0}^{-1}(T_0)$. For $k \in \mathbb{N}_0$, let $\{ {T}_0 , RT^k({T}_0) , {\Sigma}_0 \}$ be the Raviart--Thomas finite element and $I_{T_0}^{RT^k}$ the local interpolation operator defined in \eqref{RTlocal}. Let $\ell$ be such that $0 \leq \ell \leq k$. For any $\hat{v} \in W^{\ell+1,p}(\widehat{T})^3$ with ${v} = ({v}_1,v_2,{v}_3)^T := {\Psi}_{T} \hat{v}$ and ${v}_0 = (v_{0,1}, v_{0,2},v_{0,3})^T := \Psi_{T_0}{v}$, it holds that
\begin{align}
\displaystyle
&\| I_{T_0}^{RT^k} v_0 - v_0 \|_{L^p(T_0)^3} 
\leq c \frac{H_{T_0}}{h_{T_0}} \Biggl(  h_{T_0} \sum_{k=1}^3 \sum_{|\varepsilon| = \ell} \alpha^{\varepsilon} \left \| \partial_{r^{(0)}}^{\varepsilon} \frac{\partial v_0}{\partial r_k^{(0)}}  \right \|_{L^p(T_0)^3} \Biggr). \label{RT616}
\end{align}
If  Assumption \ref{ass1} is imposed, it holds that
\begin{align}
\displaystyle
&\| I_{T_0}^{RT^k} v_0 - v_0 \|_{L^p(T_0)^3} \notag \\
&\hspace{1.0cm} \leq c \frac{H_{T_0}}{h_{T_0}} \Biggl( \sum_{|\varepsilon| = \ell+ 1} \mathscr{H}^{\varepsilon} \| \partial^{\varepsilon}_{x} (\Psi_{T_0}^{-1}  v_0 ) \|_{L^p( \Phi_{T_0}^{-1} (T_0))^3} \notag \\
&\hspace{1.5cm} \quad + h_{T_0} \sum_{k=1}^3 \sum_{|\varepsilon| = \ell} \mathscr{H}^{\varepsilon} \left \| \partial^{\varepsilon}_{x} \frac{\partial (\Psi_{T_0}^{-1} v_0)}{\partial r_k} \right \|_{L^p(\Phi_{T_0}^{-1}(T_0))^3} \Biggr). \label{RT616b}
\end{align}

\end{thr}

\begin{pf*}
An analogous proof of Theorem \ref{thr1161} yields the desired results \eqref{RT616} and \eqref{RT616b}, where we use the component-wise stability \eqref{rt5} and Lemma \ref{lem1143} instead of Lemma \ref{lem1142}.
\qed	
\end{pf*}

\begin{acknowledgements}
We would like to thank the anonymous referee for the valuable comments.
%We thank the editor of the journal for his useful comments.
%We thank Stuart Jenkinson, PhD, from Edanz Group (https://jp.edanz.com/ac) for editing a draft of this manuscript. ???
%I have added this acknowledgment to comply with international publishing guidelines (e.g., ICMJE, COPE, EASE) on declaring support given to authors. Please contact us if your target journal asks for a signed letter granting my permission to be acknowledged.

%We thank Professor Norikazu Saito (Tokyo University, Japan) for the useful 	comments.%We are thankful for the time and energy you expended.
%We thank Glenn Pennycook, MSc, from Edanz Group (www.edanzediting.com/ac) for editing a draft of this manuscript.
\end{acknowledgements}




\end{document}